\documentclass[preprint2]{aastex}

\usepackage{color}
\usepackage{epsfig}
\usepackage[fleqn]{amsmath}
\usepackage{multirow}
\usepackage{multicol}
\usepackage{mathptmx}
\usepackage{graphicx}
\usepackage{bm}
\usepackage{mathrsfs}

\shorttitle{} \shortauthors{Zhang et al.}

\begin{document}


\title{Adjustment of force-gradient operator in symplectic methods\large }

\author{Lina Zhang$^{1}$, Xin Wu$^{1,2,3, \dag}$, Enwei Liang$^{1,3}$}
\affil{1. School of Physical Science and Technology, Guangxi
University, Nanning 530004, China \\ 2. School of Mathematics,
Physics and Statistics $\&$ Center of Application and Research of
Computational Physics, Shanghai University of Engineering Science,
Shanghai 201620, China \\ 3. Guangxi Key Laboratory for
Relativistic Astrophysics, Guangxi University, Nanning 530004,
China} \email{$\dag$ Corresponding Author: wuxin$\_$1134@sina.com,
xinwu@gxu.edu.cn} 

\begin{abstract}
Many force-gradient explicit symplectic integration algorithms
have been designed for the Hamiltonian $H=T
(\mathbf{p})+V(\mathbf{q})$ with kinetic energy
$T(\mathbf{p})=\mathbf{p}^2/2$ in the existing references. When
the force-gradient operator is appropriately adjusted as a new
operator, they are still suitable for a class of Hamiltonian
problems $H=K(\mathbf{p},\mathbf{q})+V(\mathbf{q})$ with
\emph{integrable} part $K(\mathbf{p},\mathbf{q}) = \sum_{i=1}^{n}
\sum_{j=1}^{n}a_{ij}p_ip_j+\sum_{i=1}^{n} b_ip_i$, where
$a_{ij}=a_{ij}(\textbf{q})$ and $b_i=b_i(\textbf{q})$ are
functions of coordinates $\textbf{q}$. The newly adjusted operator
is not a force-gradient operator but is similar to the
momentum-version operator associated to the potential $V$. The
newly extended (or adjusted) algorithms are no longer solvers of
the original Hamiltonian, but are solvers of slightly modified
Hamiltonians. They are explicit symplectic integrators with
symmetry or time-reversibility. Numerical tests show that the
standard symplectic integrators without the new operator are
generally poorer than the corresponding extended methods with the
new operator in computational accuracies and efficiencies. The
optimized methods have better accuracies than the corresponding
non-optimized counterparts. Among the tested symplectic methods,
the two extended optimized seven-stage fourth-order methods of
Omelyan, Mryglod and Folk exhibit the best numerical performance.
As a result, one of the two optimized algorithms is used to study
the orbital dynamical features of a modified H\'{e}non-Heiles
system and a spring pendulum. These extended integrators allow for
integrations in Hamiltonian problems, such as the spiral structure
in self-consistent models of rotating galaxies and the spiral arms
in galaxies.


\end{abstract}



\maketitle%

\section{Introduction}

In some cases, many trajectories of nonlinear ordinary
differential dynamical systems exhibit chaotical behavior; namely,
the separations between the trajectories and their nearby
trajectories display exponentially sensitive dependence on initial
conditions. The trajectories are not analytically integrable, and
therefore their studies mainly rely on numerical integration
schemes.  In general, detecting the chaotical behavior needs long
enough numerical integrations with reliable results. Thus, the
adopted computational schemes are demanded to perform good
stability and high precision. An eighth- and ninth-order
Runge-Kutta-Fehlberg integrator [RKF8(9)] with adaptive step sizes
significantly improves the accuracies  of the integrals of motion
and the trajectories. However, the higher-precision solutions come
at the expense of computational time. Furthermore, lower-order
numerical methods that preserve geometric properties of the flow
of differential equations, i.e., geometric numerical integration
schemes [1], are often employed to achieve very accurate long-time
determination of the trajectories. The geometric integrators
include symplectic integrators for Hamiltonian systems [2-4],
symmetric integrators for reversible systems, manifold correction
schemes for the consistency of integrals (or quasi integrals) of
motion [5, 6], energy-preserving methods [7], etc.

When dealing with Hamiltonian systems, symplectic integrators are
the most appropriate geometric solvers, which own the symplectic
nature of Hamiltonian dynamics. The errors in the integrals of
motion involving energy integral have no secular growth and tend
to zero for infinitesimal time steps. In this sense, symplectic
methods approximately conserve the integrals of motion. There are
implicit symplectic schemes for inseparable Hamiltonian systems,
and explicit symplectic schemes for integrable separable
Hamiltonian systems. The implicit schemes were developed by Feng
and Qin [8] using generating functions. Implicit symplectic
Runge-Kutta methods including the second-order implicit midpoint
rule were presented by Sanz-Serna [9]. The family of Gauss-Legndre
Runge-Kutta methods [10] are also symplectic. These implicit
methods are applied to full Hamiltonian systems. If a part of the
Hamiltonian systems is explicitly solved and another part of the
Hamiltonian systems is implicitly solved, then explicit and
implicit combined symplectic methods [11-13] are obtained and can
reduce the expense of computational time, compared with the
implicit symplectic methods for the full Hamiltonian systems.

Are symplectic integrations always implicit for nonseparable
Hamiltonians? No, they are not. In fact, explicit symplectic
integrations have been possibly available for some nonseparable
Hamiltonians [14]. Relatively recently, explicit symplectic
integrators were designed for the Hamiltonian of Schwarzschild
spacetime split into four integrable parts [15]. There are a class
of explicit extended phase-space symplectic or symplectic-like
integrations [16] for an arbitrary separable or inseparable
Hamiltonian system. Of course, explicit symplectic integrations
are less computationally expensive in general than same order
implicit symplectic integrations.

Naturally, explicit symplectic integrations are extensively
applied to separable Hamiltonians.  Ruth [3] proposed second- and
third-order explicit symplectic methods for Hamiltonian systems of
the form $H=T(\mathbf{p})+V(\mathbf{q})$. Along this direction,
higher order standard explicit symplectic schemes were developed
by many authors [17-20]. It is worth emphasizing that Ruth [3]
also gave another third-order symplectic algorithm in which the
computation of force gradient is included. Such construction is
the force-gradient symplectic integrator. Based on this idea, a
series of higher order explicit force-gradient symplectic
integration algorithms were established and applied by several
authors [21-23]. Positive time coefficients can be admissible in many of the algorithmic constructions. In view
of this, the force-gradient algorithms with
positive intermediate time-steps are suitable for solving
time-irreversible problems, such as imaginary time Schr\"{o}dinger
equations. A force-gradient explicit symplectic integration needs
less exponential functions of Lie operators than a same order
standard explicit symplectic integration. The former algorithm is
generally superior to the latter one at same order in accuracy.

When dealing with Hamiltonian systems with the integrable
perturbation decomposition form
$H=H_0(\mathbf{p},\mathbf{q})+\epsilon H_1(\mathbf{q})$, the
Wisdom-Holman symplectic map of second order [24] drastically
improves the numerical accuracy, compared with the standard
explicit symplectic method  of second order for the Hamiltonian
splitting form $H=T(\mathbf{p})+V(\mathbf{q})$. According to the
perturbation decomposition of  Hamiltonian systems, a number of
higher order explicit symplectic integrations were developed and
generalized in some references [25, 26].

It can be seen from the above presentations that the construction
of symplectic integrators is closely  related to Hamiltonian
systems or their splitting forms. In this contribution, we plan to
develop explicit force-gradient symplectic integration algorithms
for Hamiltonian systems of the form
$H=K(\mathbf{p},\mathbf{q})+V(\mathbf{q})$ with integrable kinetic
energy $K(\mathbf{p},\mathbf{q}) = \sum_{i=1}^{n}
\sum_{j=1}^{n}a_{ij}{}p_ip_j+\sum_{i=1}^{n} b_ip_i$, where
$a_{ij}=a_{ij}(\textbf{q})$ and $b_i=b_i(\textbf{q})$. Because
this Hamiltonian is different from the Hamiltonian
$H=T(\mathbf{p})+V(\mathbf{q})$,  the force-gradient operator for
the latter Hamiltonian should be modified appropriately in the
former Hamiltonian.

The remainder of the present paper is organized as follows. The
force-gradient operator is extended or adjusted in Sect. 2. Taking
two models, we estimate the numerical performance of some extended
version fourth-order force-gradient symplectic integrators in
Sect. 3. For comparison, the standard fourth-order symplectic
method of Forest $\&$ Ruth [17] and its optimized methods [20] are
employed. Possible differences in the dynamics of order and chaos
among the tested symplectic integrators are compared, and the
dynamics of the two systems is truly described. Finally, the main
results are concluded in Sect. 4.

\section{Extension of force-gradient operator}

In this Section, some known force-gradient symplectic integrators
are introduced. Then, their applications are extended. In what
follows, the symplecticity of the extended algorithms is clearly
shown.

\subsection{Existing force-gradient symplectic
integrators }

A kinetic energy $T(\mathbf{p})$ is a quadratic function of
$n$-dimensional momentum vector $\mathbf{p}=(p_1,\cdots, p_n)$
\begin{equation}
T(\mathbf{p})=\frac{1}{2} \mathbf{p}^2,
\end{equation}
and  a potential energy $V(\mathbf{q})$ depends on coordinate
$\mathbf{q}=(q_1,\cdots, q_n)$ only. $T$ and $V$ determine a
Hamiltonian
\begin{equation}
H(\mathbf{p},\mathbf{q})=T(\mathbf{p})+V(\mathbf{q}).
\end{equation}
Lie derivative operators of $T$ and $V$ are defined as
\begin{eqnarray}
  A &=& \{,T\} = \sum_{i=1}^{n}T_{p_{i}} \frac{\partial}{\partial q_{i}}, \\
  B &=& \{,V\} = -\sum_{i=1}^{n} V_{q_{i}} \frac{\partial}{\partial
  p_{i}},
\end{eqnarray}
where $T_{p_{i}}=\partial T / \partial p_{i}=p_{i}$,
$V_{q_{i}}=\partial V / \partial q_{i}$, and symbols $\{, ~\}$
denote Poisson brackets. Applying Lie derivative $A$ to act on
coordinates $q_{i}$ and momenta $p_{i}$, we obtain $\dot{q}_{i}=A
q_{i}=T_{p_{i}}=p_{i}$ and $\dot{p}_{i}=A p_{i}=0$. Clearly, $A$
is a position-version operator. In other words, $T$ has an
analytical solution as an explicit function of time. Because
$Bp_{i}=- V_{q_{i}}$ and  $Bq_{i}=0$, $B$ is a momentum-version
operator and is easily analytically solved.

The two operators $A$ and $B$ can symmetrically compose a
second-order Verlet symplectic integrator [27]
\begin{eqnarray}
  M2=e^{W}=e^{\frac{\tau}{2} B} e^{\tau A} e^{\frac{\tau}{2} B},
\end{eqnarray}
where $\tau$ is a step size, and $W$ is written in terms of
Baker-Campbell-Hausdroff (BCH) formula as
\begin{eqnarray}
 W &=& \tau(A+B)+\tau^{3}(-\frac{1}{12}[A,[B, A]] \nonumber \\
 && +\frac{1}{24}[B,[A, B]])+\mathcal{O}\left(\tau^{5}\right).
\end{eqnarray}
The two operators can also symmetrically compose fourth-order
explicit symplectic algorithms, such as the Forest-Ruth method
[17]
\begin{eqnarray}
 M4 &=& \nonumber e^{\alpha \tau A} e^{\beta \tau B} e^{(\frac{1}{2}-\alpha) \tau A} e^{(1-2 \beta) \tau B}\\
&&  \times e^{(\frac{1}{2}-\alpha) \tau A} e^{\beta \tau B}
e^{\alpha \tau
 A},
\end{eqnarray}
where $\beta=1 /(2-\sqrt[3]{2})$ and $\alpha=\beta/2$. An
optimized Forest-Ruth-like explicit symplectic  algorithm of order
4 in [20] is
\begin{eqnarray}
M4V &=& \nonumber e^{\xi \tau B}e^{(1-2 \lambda) \frac{\tau}{2} A}
e^{\chi \tau B} e^{\lambda \tau A} e^{(1-2(\chi+\xi)) \tau B} \\
&& \times e^{\lambda \tau A}e^{\chi \tau B}e^{(1-2 \lambda)
\frac{\tau}{2} A}e^{\xi \tau B},
\end{eqnarray}
where time coefficients are
\begin{eqnarray}
\xi &=& 0.1644986515575760E+00, \nonumber \\
\lambda &=& -0.2094333910398989E-01, \nonumber \\
\chi &=& +0.1235692651138917E+01. \nonumber
\end{eqnarray}
In fact, only two of the three coefficients sufficiently satisfy
the conditions for order 4. In this case, one of the three
coefficients has a free choice. However, the free coefficient such
as $\xi$ can be determined when the norm of the leading term of
fifth-order truncation errors is minimized. This is the so-called
optimized method. Another optimized algorithm of order 4 is
\begin{eqnarray}
M4P &=& \nonumber e^{\xi \tau A}e^{(1-2 \lambda) \frac{\tau}{2} B}
e^{\chi \tau A} e^{\lambda \tau B} e^{(1-2(\chi+\xi)) \tau A} \\
&& \times e^{\lambda \tau B}e^{\chi \tau A}e^{(1-2 \lambda)
\frac{\tau}{2} B}e^{\xi \tau A},
\end{eqnarray}
where time coefficients are
\begin{eqnarray}
\xi &=& 0.1786178958448091E+00, \nonumber \\
\lambda &=& -0.2123418310626054E+00, \nonumber \\
\chi &=& -0.6626458266981849E-01. \nonumber
\end{eqnarray}
Explicit symplectic algorithms to arbitrary even orders were given
in [18, 19].

Let us consider one of the truncation error terms about $\tau^{3}$
in Eq. (6). In terms of commutator $[A,B]=AB-BA$, we have a
commutator $C=[B,[A, B]]=[B, AB-BA]=2BAB-BBA-ABB$. It is easy to
derive the following results: $ABBq_{i}=ABBp_{i} \equiv 0$,
$BABq_{i} \equiv 0$, $BBAq_{i}=BBAp_{i} \equiv 0$, and
\begin{eqnarray}
  BABp_{i} = \sum_{j=1}^{n} V_{q_{i}
q_{j}} V_{q_{j}}.
\end{eqnarray}
Thus, $C$ is still a momentum-version operator
\begin{eqnarray}
  C &=& \nonumber 2BAB \\ \nonumber
  &=& \sum_{i=1}^{n} \sum_{j=1}^{n}2 V_{q_i q_j} V_{q_j}\frac{\partial}{\partial
  p_{i}} \nonumber \\
  &=& \sum_{i=1}^{n} \nabla_i \mathbf{f}^2 \frac{\partial}{\partial
  p_{i}},
\end{eqnarray}
where $V_{q_i q_j}=\partial^{2}V/\partial q_i \partial q_j$, and
$\mathbf{f}=(f_1,\cdots, f_n)$ with $\mathbf{f}_{j}=-V_{q_{j}}$ is
a component of the force governed by potential $V$. That is to
say, $C$ is a momentum-version operator with respect to the
gradient of the square of force. In this sense, $C$ is called as a
force-gradient operator in Refs. [21-23].

The above demonstrations show that the third-order truncation
error term $C$ in the second-order method M2 and $B$ belong to
momentum-version operators. This means that a combination of $B$
and $C$ can yield higher-order algorithms.  For example,
\begin{eqnarray}
 F2=e^{\frac{\tau}{2} (B+\frac{\tau^{2}}{24}C)} e^{\tau A} e^{\frac{\tau}{2} (B+\frac{\tau^{2}}{24}C)}
\end{eqnarray}
is a second-order algorithm which eliminates the error term
$[B,[A,B]]$ in Eq. (6). If the error term $[A,[B,A]]$ in Eq. (6)
is also eliminated, a five-stage fourth-order method is Chin's
construction [21]
\begin{eqnarray}
F4=\nonumber e^{\frac{\tau}{2}\left(1-\frac{1}{\sqrt{3}}\right) A} e^{\frac{\tau}{2} (B+\frac{\tau^{2}}{24}(2-\sqrt{3}) C)} \\
\times e^{\frac{\tau}{\sqrt{3}}A} e^{\frac{\tau}{2}
(B+\frac{\tau^{2}}{24}(2-\sqrt{3}) C)}
e^{\frac{\tau}{2}\left(1-\frac{1}{\sqrt{3}}\right) A}.
\end{eqnarray}
When $A\leftrightarrow B$ in Eq. (13), there is another five-stage
fourth-order scheme
\begin{eqnarray}
F4^{*}= e^{\frac{\tau}{6} (B+\frac{\tau^2}{72}C)}e^{\frac{\tau}{2}
A}e^{\frac{2\tau}{3} [B+\frac{\tau^2}{72}C]} e^{\frac{\tau}{2}
A}e^{\frac{\tau}{6} [B+\frac{\tau^2}{72}C]}.
\end{eqnarray}
An optimized five-stage fourth-order method with the operator $C$
[23] is
\begin{eqnarray}
F4O &=& e^{\lambda\tau [B+(6\xi+k)\tau^2C]}e^{\frac{\tau}{2}
A}e^{(1-2\lambda)\tau [B+(6\xi+k)\tau^2C]} \nonumber \\ && \times
e^{\frac{\tau}{2} A}e^{\lambda\tau [B+(6\xi+k)\tau^2C]},
\end{eqnarray}
where $\lambda=1/6$, $\xi=-17/18000$, $\chi=71/4500$, and
$k=-4\xi-\frac{1}{2}+\frac{3\chi}{2}$. The choice of $\xi$
minimizes the norm of the leading term of fifth-order truncation
errors. Two optimized seven-stage fourth-order methods with the
operator $C$ in [23] are
\begin{eqnarray}
F4V &=& e^{\lambda\tau [B+2(\xi+\chi)\tau^2C]}e^{\theta\tau
A}e^{(1-2\lambda)\frac{\tau}{2} [B+2(\xi+\chi)\tau^2C]} \nonumber
\\ && \times e^{(1-2\theta)\tau A} e^{(1-2\lambda)\frac{\tau}{2}
[B+2(\xi+\chi)\tau^2C]}e^{\theta\tau A} \nonumber
\\ && \times e^{\lambda\tau
[B+2(\xi+\chi)\tau^2C]},
\end{eqnarray}
where
\begin{eqnarray}
\theta &=& 0.2728983001988755E+00, \nonumber \\
\lambda &=& 0.8002565306418866E-01, \nonumber \\
\chi &=& 0.2960781208329478E-02, \nonumber \\
\xi &=& 0.2725753410753895E-03, \nonumber
\end{eqnarray}
and
\begin{eqnarray}
F4P &=& e^{\theta\tau A}e^{\lambda\tau [B+(2\xi+\chi)\tau^2C]}
e^{(1-2\theta)\frac{\tau}{2} A} \nonumber \\ && \times
e^{(1-2\lambda)\tau [B+(2\xi+\chi)\tau^2C]}
e^{(1-2\theta)\frac{\tau}{2} A} \nonumber \\ && \times
e^{\lambda\tau [B+(2\xi+\chi)\tau^2C]}e^{\theta\tau A},
\end{eqnarray}
where
\begin{eqnarray}
\theta &=& 0.1159953608486416E+00, \nonumber \\
\lambda &=& 0.2825633404177051E+00, \nonumber \\
\chi &=& 0.3035236056708454E-02, \nonumber \\
\xi &=& 0.1226088989536361E-02. \nonumber
\end{eqnarray}
Note that the combinations of $B$ and $C$ in Eqs. (15)-(17) are
slightly unlike those in Ref. [23], but the time coefficients are
the same as those in Ref. [23]. Although the three methods F4O,
F4V and F4P are optimized algorithms, their differences are that a
free coefficient for F4O minimizes the norm of the leading term of
fifth-order truncation errors, but two free coefficients such as
$\xi$ and $\chi$ for F4V and F4P do.

Algorithms (12)-(17) with the force-gradient operator $C$ are the
so-called force-gradient explicit symplectic methods in the known
publications (e.g., Refs. [21-23]). Their constructions are based
on  one of the third-order truncation error terms in the
second-order method M2 included in the potential energy.
Only the time coefficient $\lambda(6\xi+k)$ with
respect to the operator $C$ in Eq. (15) is negative, but the time
coefficients with respect to the operators $A$ and $B$ in Eqs.
(12)-(17) and the time coefficients with respect to the operator
$C$ in Eqs. (12)-(14), (16) and (17) are positive. Thus, no F4O
but F2, F4, F4$^*$, F4V and F4P can have positive intermediate
time-steps.

\subsection{Adjustment of the force-gradient operator}

Now, let us suppose a Hamiltonian
\begin{eqnarray}
J(\mathbf{p},\mathbf{q}) &=&
K(\mathbf{p},\mathbf{q})+V(\mathbf{q}), \\
 K(\mathbf{p},\mathbf{q}) &=& \sum_{i=1}^{n}
\sum_{j=1}^{n}a_{ij}p_ip_j+\sum_{i=1}^{n} b_ip_i,
\end{eqnarray}
where $a_{ij}=a_{ij}(\textbf{q})$ and $b_i=b_i(\textbf{q})$.
$K(\mathbf{p},\mathbf{q})$ is a formal kinetic energy unlike the
kinetic energy in Eq. (1). It is a quadratic function of the
momentum vector $\mathbf{p}$, and also depends on the coordinate
vector $\mathbf{q}$.

Lie derivative operator of $K$ is expressed as
\begin{eqnarray}
  A &=& \{,K\} = \sum_{i=1}^{n}\left(K_{p_{i}} \frac{\partial}{\partial q_{i}}
  -K_{q_{i}} \frac{\partial}{\partial p_{i}}\right),
\end{eqnarray}
where $K_{p_{i}}=\partial K / \partial p_{i}$ and
$K_{q_{i}}=\partial K / \partial q_{i}$. We have $\dot{q}_{i}=A
q_{i}=K_{p_{i}}$ and $\dot{p}_{i}=A p_{i}=-K_{q_{i}}$. In this
sense, $A$ is a momentum- and position-version mixed operator. If
$K$ is nonlinear, it is not necessarily analytically solved. Here,
$K$ is assumed to have an analytical solution as an explicit
function of time, i.e., operator $A$ is analytically solvable.
Taking $D=[B,[A, B]]=[B, AB-BA]=2BAB-BBA-ABB$, we obtain the
results as follows: $ABBq_{i}=ABBp_{i} \equiv 0$, $BABq_{i} \equiv
0$, $BBAq_{i} \equiv 0$, and
\begin{eqnarray}
  BABp_{i} &=& \sum_{j=1}^{n} \sum_{k=1}^{n} V_{q_{i}
q_{j}} V_{q_{k}} K_{p_{j} p_{k}}, \\
BBAp_{i} &=& -\sum_{j=1}^{n} \sum_{k=1}^{n} V_{q_{j}} V_{q_{k}}
K_{q_{i} p_{j} p_{k}},
\end{eqnarray}
where $K_{p_{j} p_{k}}=\frac{\partial^{2}K}{\partial p_j
\partial p_k}$ and $K_{q_{i} p_{j} p_{k}}=\frac{\partial^{3}K}{\partial
q_i\partial p_j\partial p_k}$. Thus, $D$ is still a
momentum-version operator
\begin{eqnarray}
  D &=& \nonumber 2BAB-BBA \\ \nonumber
  &=& \sum_{i=1}^{n} \sum_{j=1}^{n} \sum_{k=1}^{n}(2 V_{q_i q_j} V_{q_k}
  K_{p_jp_k}  \nonumber \\
  && +V_{q_j} V_{q_k} K_{q_i p_j p_k}) \frac{\partial}{\partial
  p_{i}}.
\end{eqnarray}
Obviously, $D$ is different from the force-gradient operator $C$.
It is an extended or adjusted operator of the force-gradient
operator $C$.

Algorithms (12)-(17) with the force-gradient operator $C$ become
useless for the Hamiltonian $J$ without doubt. However, when the
force-gradient operator $C$ is replaced with the extended operator
$D$, these integrators are still available  for the Hamiltonian
$J$. That is, they correspond in sequence to the following forms
\begin{eqnarray}
 N2=e^{\frac{\tau}{2} (B+\frac{\tau^{2}}{24}D)} e^{\tau A} e^{\frac{\tau}{2}
 (B+\frac{\tau^{2}}{24}D)},
\end{eqnarray}
\begin{eqnarray}
N4=\nonumber e^{\frac{\tau}{2}\left(1-\frac{1}{\sqrt{3}}\right) A}
e^{\frac{\tau}{2} (B+\frac{\tau^{2}}{24}(2-\sqrt{3}) D)} \\
\times e^{\frac{\tau}{\sqrt{3}}A} e^{\frac{\tau}{2}
(B+\frac{\tau^{2}}{24}(2-\sqrt{3}) D)}
e^{\frac{\tau}{2}\left(1-\frac{1}{\sqrt{3}}\right) A},
\end{eqnarray}
\begin{eqnarray}
N4^{*}= e^{\frac{\tau}{6} (B+\frac{\tau^2}{72}D)}e^{\frac{\tau}{2}
A}e^{\frac{2\tau}{3} [B+\frac{\tau^2}{72}D]} e^{\frac{\tau}{2}
A}e^{\frac{\tau}{6} [B+\frac{\tau^2}{72}D]},
\end{eqnarray}
\begin{eqnarray}
N4O &=& e^{\lambda\tau [B+(6\xi+k)\tau^2D]}e^{\frac{\tau}{2}
A}e^{(1-2\lambda)\tau [B+(6\xi+k)\tau^2D]} \nonumber \\ && \times
e^{\frac{\tau}{2} A}e^{\lambda\tau [B+(6\xi+k)\tau^2D]},
\end{eqnarray}
\begin{eqnarray}
N4V &=& e^{\lambda\tau [B+2(\xi+\chi)\tau^2D]}e^{\theta\tau
A}e^{(1-2\lambda)\frac{\tau}{2} [B+2(\xi+\chi)\tau^2D]} \nonumber
\\ && \times e^{(1-2\theta)\tau A} e^{(1-2\lambda)\frac{\tau}{2}
[B+2(\xi+\chi)\tau^2D]}e^{\theta\tau A} \nonumber
\\ && \times e^{\lambda\tau
[B+2(\xi+\chi)\tau^2D]},
\end{eqnarray}
\begin{eqnarray}
N4P &=& e^{\theta\tau A}e^{\lambda\tau [B+(2\xi+\chi)\tau^2D]}
e^{(1-2\theta)\frac{\tau}{2} A} \nonumber \\ && \times
e^{(1-2\lambda)\tau [B+(2\xi+\chi)\tau^2D]}
e^{(1-2\theta)\frac{\tau}{2} A} \nonumber \\ && \times
e^{\lambda\tau [B+(2\xi+\chi)\tau^2D]}e^{\theta\tau A}.
\end{eqnarray}
Higher-order force-gradient algorithms in Refs. [22, 23] can
become the extended methods similar to the constructions
(25)-(29).

The time coefficients of the composite operators $A$, $B$ and $D$
in each of the newly extended methods (24)-(29) correspond to
those in the existing force-gradient algorithms (12)-(17).
No N4O but N2, N4, N4$^*$, N4V and N4P can have
positive intermediate time-steps. These extended methods are
symmetric or time-reversible. The condition for symmetry or
time-reversibility is, e.g., $N4O(\tau)\times N4O(-\tau)=id$ for
the method N4O [1], where \emph{id} denotes an identical map. The
condition is easily checked from a theoretical point of view
because the exponents of these operators are linear combinations
of $\tau$ and $\tau^3$ terms. If the exponents include even power
terms of $\tau$ such as $\tau^2$ and $\tau^4$ terms, then the
condition for symmetry or time-reversibility is not satisfied.

\subsection{Preservations of symplecticity and volume of the phase space}

The standard algorithms M2, M4, M4V and M4P are
symplectic because $A$ as a Lie operator with respect to the phase
flow of a sub-Hamiltonian $H_1=T(\mathbf{p})$ is symplectic, and
$B$ as a Lie operator with respect to the phase flow of another
sub-Hamiltonian $H_2=V(\mathbf{q})$ is also symplectic. A product
of the two symplectic operators and its compositions are still
symplectic for the Hamiltonian (2). That is to say, when operators
$A$ and $B$ correspond to the phase flows of the two
sub-Hamiltonians for the Hamiltonian (2) (namely, the Hamiltonian
(2) is a symplectically separable Hamiltonian system), their
composition products are naturally symplectic from a physical
point of view [28]. Noting this idea, we can similarly show the
symplecticity of the existing force-gradient algorithms (12)-(17)
and the extended algorithms (24)-(29). Some discussions are given
as follows.

At first, we investigate the construction mechanisms of the
existing force-gradient algorithms (12)-(17). As is mentioned
above, they are constructed by adding one of the third-order
truncation error terms in the second-order method M2 to the
potential energy. This is one path for understanding the
construction mechanisms of the original force-gradient algorithms.
Another path is that the force-gradient algorithms are not
directly applied to solve the Hamiltonian (18) but are applied to
solve some modified Hamiltonians. To show this, we take a
Hamiltonian solved by the algorithm F2 as an example. Eq. (12) is
a solver of a modified Hamiltonian $H_{F2}=T+V_{F2}$, where a
modified potential is $V_{F2}=V(\mathbf{q})-\frac{\tau^2}{24}
\mathbf{f}^2$. Operator $A$ acting on the sub-Hamiltonian $T$
advancing a time of $\tau$ is the exact phase flow $e^{\tau A}$ of
the sub-Hamiltonian $T$, and is symplectic. Operator
$B_{F2}=(B+\frac{\tau^{2}}{24}C)$ acting on the sub-Hamiltonian
$V_{F2}$ advancing a time of $\tau/2$ is the exact phase flow
$e^{\frac{\tau}{2}B_{F2}}$ of the sub-Hamiltonian $V_{F2}$, and is
also symplectic. Naturally, a composition porduct of these
symplectic operators such as $e^{\tau A}e^{\frac{\tau}{2}B_{F2}}$
is symplectic [28]. In fact, the explicit symmetric composition
scheme $F2$ for the symplectically separable modified Hamiltonian
$H_{F2}$  is the standard second-order symplectic method M2 with
$B\rightarrow B_{F2}$. Of course, different force-gradient
symplectic algorithms such as F4 correspond to different modified
Hamiltonians. Namely, the force-gradient symplectic algorithms for
the original Hamiltonian (2) are the standard symplectic methods
for the corresponding modified Hamiltonians.

The construction mechanisms of the extended algorithms (24)-(29)
are similar to those of the existing force-gradient algorithms
(12)-(17). To show this, we consider the algorithm N2 for the
Hamiltonian problem $J$ in Eq. (18). There is a modified potential
$V_{N2}(\mathbf{q})=V(\mathbf{q})-\frac{\tau^2}{24}V_D(\mathbf{q})$,
where $V_D(\mathbf{q})$ is solved by $\partial
V_D(\mathbf{q})/\partial q_i= \sum_{j=1}^{n} \sum_{k=1}^{n}(2
a_{jk} V_{q_i q_j} V_{q_k}
  +V_{q_j} V_{q_k} \frac{\partial a_{jk}}{\partial q_i})$.
Although the function $V_D$ is difficulty written in general, it
should exist from a mathematical point of view. In fact, we do not
need to know what the function $V_D$ is, but we need to know only
what the function $\partial V_D/\partial q_i$ is. The modified
Hamiltonian solved by the algorithm N2 is a
symplectically separable Hamiltonian system
$J_{N2}=K+V_{N2}(\mathbf{q})$. Only one difference in the
construction mechanisms between the force-gradient algorithm F2
and the extended algorithm N2 is that the modified potential
$V_{F2}$ for F2 is easily expressed, whereas the modified
potential $V_{N2}$ for N2 is not. This does not destroy the
symplecticity of the algorithm N2 with the extended operator $D$
for the Hamiltonian (18). Similarly, the modified Hamiltonians
solved by the algorithms N4 and $N4^{*}$ can be expressed as
symplectically separable Hamiltonian systems
$J_{N4}=K+[V(\mathbf{q})-\frac{\tau^2}{24}(2-\sqrt{3})V_D(\mathbf{q})]$
and $J_{N4^*}=K+[V(\mathbf{q})-\frac{\tau^2}{72}V_D(\mathbf{q})]$.
In this way, the symplecticity of the algorithms N4 and $N4^{*}$
is shown sufficiently. The modified Hamiltonians for the extended
algorithms (27)-(29) also exist and therefore the extended
algorithms  remain symplectic.

Precisely speaking, the  symplecticity of each of
the aforementioned algorithms should satisfy the condition
\begin{eqnarray}
\textbf{S}^{T} \textbf{I} \textbf{S}=\textbf{I},
\end{eqnarray}
where $\textbf{S}$ and $\textbf{I}$ are $2n\times 2n$ matrixes.
Here, we take four-dimensional phase-space variables
$\textbf{Z}=(x, y, p_x, p_y)$ as an example to show the
expressions of $\textbf{S}$ and $\textbf{I}$. The solutions
$\textbf{Z}_m=(x_m, p_{xm}, y_m, p_{ym})$ at an $m$th step from
the solutions $\textbf{Z}_{m-1}=(x_{m-1}, p_{xm-1}$, $y_{m-1},
p_{ym-1})$ at an $(m-1)$th step advancing the time step $\tau$ for
one of the aforementioned algorithms  are expressed as
\begin{eqnarray}
\textbf{Z}^{T}_m=\textbf{f}^{T}(\tau, \textbf{Z}_{m-1}),
\end{eqnarray}
where $\textbf{f}=(f_x, f_{px},f_y, f_{py})$; that is,
\begin{eqnarray}
x_m &=& f_x (\tau, x_{m-1}, p_{xm-1}, y_{m-1}, p_{ym-1}), \nonumber \\
p_{xm} &=& f_{px} (\tau, x_{m-1}, p_{xm-1}, y_{m-1}, p_{ym-1}), \nonumber \\
y_m &=& f_y (\tau, x_{m-1}, p_{xm-1}, y_{m-1}, p_{ym-1}), \nonumber \\
p_{xm} &=& f_{py} (\tau, x_{m-1}, p_{xm-1}, y_{m-1}, p_{ym-1}).
\nonumber
\end{eqnarray}
Their differential forms are
$d\textbf{Z}^{T}_m=\textbf{S}d\textbf{Z}^{T}_{m-1}$, where
$\textbf{S}$ is a $4\times 4$ matrix
\begin{equation}
 \textbf{S}= \left(\begin{array}{cccc}
 $$ \frac{\partial x_{m}}{\partial x_{m-1}}&\frac{\partial x_{m}}{\partial y_{m-1}}&\frac{\partial x_{m}}{\partial p_{xm-1}}&\frac{\partial x_{m}}{\partial p_{ym-1}}\\
   \frac{\partial p_{xm}}{\partial x_{m-1}}&\frac{\partial p_{xm}}{\partial y_{m-1}}&\frac{\partial p_{xm}}{\partial p_{xm-1}}&\frac{\partial p_{xm}}{\partial p_{ym-1}}\\
   \frac{\partial y_{m}}{\partial x_{m-1}}&\frac{\partial y_{m}}{\partial y_{m-1}}&\frac{\partial y_{m}}{\partial p_{xm-1}}&\frac{\partial y_{m}}{\partial p_{ym-1}}\\
  \frac{\partial p_{ym}}{\partial x_{m-1}}&\frac{\partial p_{ym}}{\partial y_{m-1}}&\frac{\partial p_{ym}}{\partial p_{xm-1}}&\frac{\partial p_{ym}}{\partial p_{ym-1}}$$
  \end{array}\right).
\end{equation}
$\textbf{I}$ is another $4\times 4$ matrix
\begin{equation}
 \textbf{I}= \left(\begin{array}{cccc}
 $$ 0 &  0 & 0 & -1 \\
   0 & 0& -1 & 0 \\
   0 &  1 & 0 & 0 \\
   1 & 0& 0 & 0 $$
  \end{array}\right).
\end{equation}
Matrix $\textbf{S}$ for satisfying the condition (30) is a
symplectic matrix. Such symplectic matrixes are present for the
above algorithms M2, $\cdots$, F2, $\cdots$, N2,  $\cdots$. This
symplecticity means a symplectic structure described by a closed
nondegenerate differential 2-form $\omega=dx\wedge dp_{x} +
dy\wedge dp_{y}$ on a four-dimensional differential manifold
$M_{4}$ (symbol $\wedge$ denotes a wedge product). The condition
(30) corresponds to the equality of the symplectic structure
$\omega_{m-1}$ at the $(m-1)$th step and the symplectic structure
$\omega_{m}$ at the $m$th step: $\omega_m=\omega_{m-1}$. A
differential 4-form comes from the product of original
differential 2-forms: $\widetilde{\omega}=dx\wedge dp_{x}\wedge
dy\wedge dp_{y}$. The preservation of differential 2-form
naturally leads to that of differential 4-form
\begin{eqnarray}
\widetilde{\omega}_{m}= det|\textbf{S}|\widetilde{\omega}_{m-1},
\end{eqnarray}
where the determinant of the matrix $\textbf{S}$ is
$det|\textbf{S}|=1$. From a geometric point of view, the
differential 4-form corresponds to a volume of the phase space:
\begin{eqnarray}
&& Vol=|\widetilde{\omega}_{m}|= |det|\textbf{S}||\cdot
|\widetilde{\omega}_{m-1}|.
\end{eqnarray}
Therefore, the differential 4-form and volume of the phase space
are preserved when $det|\textbf{S}|=1$ (see [28] for more details
on the symplectic structure and volume of the phase space). The
result on $det|\textbf{S}|=1$ will be tested in later numerical
experiments.

In short, the novel contribution of this paper is to extend the
existing force-gradient symplectic algorithms for the Hamiltonian
(2) in Refs. [21-23] to solve the Hamiltonian (18). Above all, the
force-gradient operator  $C$ must be replaced by the new operator
$D$.

\section{Numerical simulations}

A modified H\'{e}non-Heiles system and a spring pendulum are taken
as two models to check the numerical performance of the extended
algorithms N4, N4O, N4V and N4P in accuracies of energy and
position. Methods M4, M4V and M4P are compared with the extended
algorithms. Possible regular and chaotic dynamical differences
between the algorithms N4 and M4 are shown. Dynamical behavior of
order and chaos in the two problems are described in terms of the
best extended algorithm.

\subsection{Modified H\'{e}non-Heiles system}

Let us consider a modified H\'{e}non-Heiles system. Set $V$ as the
potential of H\'{e}non-Heiles system [29]
\begin{eqnarray}
  V=\frac{1}{2}\left(x^{2}+y^{2}\right)+x^{2} y-\frac{1}{3}y^{3}.
\end{eqnarray}
The standard kinetic energy of H\'{e}non-Heiles system is
$T=\left(p_{x}^{2}+p_{y}^{2}\right)/2$. Here, it is slightly
modified as
\begin{eqnarray}
  K=\frac{1}{2}\left(y p_{x}^{2}+p_{y}^{2}\right).
\end{eqnarray}
$K$  belongs to one of the forms given in Eq. (19). Obviously, the
operator $A$ for $K$, the operator $B$ for $V$ and the operator
$D$ are easily, analytically solvable. Thus, the algorithms N4,
$N4^{*}$, N4O, N4V, N4P, M4, M4V and M4P are easily available for
the system $H=K+V$ (for convenience, $J$ in Eq. (18) is replaced
by $H$).

The time step is $\tau=0.1$. Energy in the system is $E=1/120$.
Initial conditions  are $x=0$, $y=-2.02$ and $p_{y}=0$. The
initial positive value of $p_{x}$ is determined by $E=H$. Fig.
1(a) plots energy errors $|\Delta H|=E_t-E$ for several
algorithms, where $E_t$ denotes the numerical energy at time $t$.
The errors  have secular growth for the conventional fourth-order
Runge-Kutta (RK4), whereas do not have for the algorithms N4, N4O,
N4V, N4P, M4, M4V and M4P. The property without secular drift in
the energy errors is due to the symplecticity of these methods.
The symplectic methods according to the energy errors listed in
Table 1 are mainly divided into three groups as follows. M4 has
the poorest anergy accuracy, and N4P and N4V have the best anergy
accuracies. N4P and N4V have almost the same accuracy, and their
accuracies are about three orders of magnitude better than the
accuracy for M4. The four methods N4, M4P, M4V and N4O have minor
differences in the energy accuracies. The errors for large to
small are N4, M4P, M4V and N4O. In particular, the accuracy for N4
is one order of magnitude better than that for M4. A position
error at time $t$ for each of the methods RK4, N4, N4O, N4V, N4P,
M4, M4V and M4P is estimated by $|\Delta
\mathbf{r}|=\sqrt{(x_2-x_1)^2+(y_2-y_1)^2}$, where the solutions
$(x_1,y_1)$ are given by the method, and the solutions $(x_2,y_2)$
are obtained from the higher-precision integrator RKF8(9). When
the integration time reaches $10^{4}$ corresponding to $10^{5}$
steps, the position errors are shown in Fig. 1(b) and Table 2. RK4
has the largest error with an order of $10^{0.47}$. The position
error for M4 is also larger than 1. The position errors for N4P
and N4V are two orders of magnitude smaller than that for  M4, and
one order of magnitude smaller than those for the four methods N4,
M4P, M4V and N4O. The four methods N4, M4P, M4V and N4O are almost
the same in the position errors. These results in Tables 1 and 2
indicate that the standard symplectic integrators without the
operator $D$ are poorer than the corresponding extended methods
with  the operator $D$ in energy accuracies (e.g., M4 inferior to
N4, M4V and M4P inferior to N4O, N4V and N4P). The optimized
methods are better than the corresponding non-optimized methods
(e.g., M4V and M4P superior to M4, and N4O, N4V and N4P superior
to N4). The optimized methods with one free coefficient are
inferior to those with two free coefficients (e.g., N4O inferior
to N4V and N4P). Clearly, M4 performs the poorest accuracy, and
N4V and N4P exhibit the best accuracies.

The smaller the time step gets, the higher the accuracy of an
integrator is. When the time step is appropriately smaller, e.g.,
$\tau=0.01$, the energy errors are described in Fig. 1(c) and
Table 1. The errors are $10^{-6.75}$ for M4, $10^{-8}$ for N4,
M4V, M4P, and N4O, and $10^{-9.7}$ for N4V and N4P. That is, the
energy errors typically decrease for each symplectic algorithm.
Although RKF8(9) as a non-symplectic method has a secular drift in
the energy errors, it has the best energy accuracy and is a good
reference integrator for evaluating the performance of the other
methods. The position errors in Fig. 1(d) and Table 2 are about
$10^{-2}$ for M4 and RK4, $10^{-5}\sim10^{-4}$ for M4V, M4P, N4
and N4O, and $10^{-6}$ for N4V and N4P. That is to say, as far as
the accuracies of energy and position for the smaller time step
$\tau=0.01$ are concerned, M4 is still the poorest one, and N4V
and N4P are the best ones among the symplectic methods. Table 3
lists CPU times for each algorithm in Figs. 1 (a) and (c). RK4 is
the fastest, and RKF8(9) is the slowest. The computational
efficiencies of the symplectic integrators from high to low are
N4, N4O, M4, N4P, N4V, M4V and M4P. Of course, there are no
relatively dramatic differences in the computational cost between
N4 (or N4O) and N4P (or N4V).

Fig. 2 is used to check whether the determinants
of the matrix $\textbf{S}$, $det|\textbf{S}|$, are 1 for the three
methods RK4, M2 and N2. The error $det|\textbf{S}|-1$ for RK4
grows with time spanning 300. This means that RK4 does not satisfy
Eqs. (34) and (35). However, the errors $det|\textbf{S}|-1$ for M2
and N2 almost remain at the machine precision in the
double-precision environment. This sufficiently supports the
preservations of the differential 4-form $\widetilde{\omega}$ and
the phase-space volume in the standard symplectic method M2 and
the extended algorithm N2. In principle, $det|\textbf{S}|\neq 1$
for  RKF8(9) and $det|\textbf{S}|= 1$ for the fourth-order methods
M4, M4V, M4P, N4, N4$^*$, N4O, N4V and N4P can be confirmed
numerically. However, the matrixes $\textbf{S}$ for these
algorithms have such long expressions (with over 1000 pages
outputted by Matlab) that their determinants are difficultly
computed.

Fig. 3 shows that the two symplectic integrators N4 and M4 provide
different dynamical phase-space structures to the same orbit in
Fig. 1 because they have different numerical accuracies for the
time step $\tau=0.1$. A single Kolmogorov-Arnold-Moser (KAM) torus
on the Poincar\'{e} section $x=0$ with $p_x>0$ is described by the
Forest-Ruth method M4 in Fig. 3(a), but many islands are given by
the newly extended method N4 in Fig. 3(b). Although the two tori
are regular, they are different. In fact, a many-islands torus
becomes easier for the occurrence of chaos than a single torus.
Which of the algorithms M4 and N4 can provide correct results? N4
can because the results of N4 are consistent with those of
higher-precision method RKF8(9) in Fig. 3(c), and are also the
same as those of the five methods M4V, M4P, N4O, N4V and N4P. The
different KAM tori described by the two methods N4 and M4 are
because N4 is one order of magnitude better than M4 in the
accuracies of energy and position, as shown in Figs. 1 (a) and
(b). In other words, the time step $\tau=0.1$ is too large to be
chosen for M4. However, the phase-space structures can be truly
described by M4 for the smaller time step $\tau=0.01$ because M4
has better accuracies for the smaller time step $\tau=0.01$ than
for the larger time step $\tau=0.1$ (see Fig. 1 and Tables 1 and 2
for more information).

Only when the initial value $y=-1.108$ with the initial value
$p_{x}$ is altered, does M4 with the larger time step $\tau=0.1$
indicate that the orbit considered in Fig. 4(a) seems to be a
many-islands torus, but exhibits the chaoticity because many
discrete points are randomly filled with small regions. Unlike M4,
N4 with the larger time step $\tau=0.1$ and RKF8(9) show the
regularity of the same orbit in Figs. 4 (b) and (c). As claimed
above, N4 is superior to M4 in accuracy, therefore, the KAM torus
is physically given by N4, but the non-physically spurious
chaoticity is caused by M4. These results are also confirmed by
fast Lyapunov indicators (FLIs) in Fig. 4(d). The FLIs are from a
modified form of Lyapunov exponents [30]. They are originally
defined in terms of the lengths of tangent vectors by
Froeschl\'{e} $\&$ Lega [31]. Using the phase-space distances
between two adjacent orbits at times 0 and $t$,  $d(0)$ and
$d(t)$, Wu et al. [32] suggested the computation of FLI according
to the following form
\begin{equation}\label{28} \emph{FLI}=\log _{10}
\frac{d(t)}{d(0)}.
\end{equation}
A bounded orbit is ordered if its FLI grows algebraically with
time $\log _{10} t$, but chaotic when its FLI increases
exponentially. In other words, the method of complete different
growth rates of FLIs with time is faster to distinguish between
the two cases of order and chaos than the technique of Lyapunov
exponents. When the integration time reaches 3000 in Fig. 4(d),
the FLI is 25 for M4, and smaller than 2.5 for N4. Clearly, M4 and
N4 for the larger time step $\tau=0.1$ indeed give the orbit
chaotic and regular dynamical behaviors, respectively. Of course,
the chaoticity for M4 is spurious because of M4 performing the
poor accuracy. If the smaller time step $\tau=0.01$ is adopted, M4
like N4 can give the true dynamical behavior to the orbit.

Given the initial value $y=-1.654$ and the larger time step
$\tau=0.1$, the methods of Poincar\'{e} section and FLI in Fig. 5
show that the integrated orbit is a regular single torus for M4,
whereas a figure-eight orbit for N4 and RKF8(9). The figure-eight
orbit seems to be regular, but is in fact chaotic due to the
existence of a hyperbolic point, which has a stable direction and
a unstable direction. The chaoticity for N4 is physical, but the
regular dynamical information given by M4 is not physical. This is
because M4 does not give high enough accuracy to the numerical
solutions, but N4 does. M4 can also provide the reliable results
for the smaller time step $\tau=0.01$.  In addition, the energy
error (not plotted) of the chaotic orbit for N4 in Fig. 5 is
similar to that of the regular orbit for N4 in Fig. 1. Namely, the
energy accuracy for N4 is independent of the regularity or
chaoticity of orbits.

Several main results can be concluded from the above
demonstrations. The standard symplectic integrators without the
operator $D$ are inferior to the corresponding extended methods
with the operator $D$ in computational accuracies and
efficiencies. The optimized methods have better accuracies than
the corresponding non-optimized methods. N4V and N4P exhibit the
best accuracies. Although N4 can provide reliable results on the
orbital dynamical behavior for the larger time step $\tau=0.1$,
one of the two optimized extended methods  N4P and N4V should be
the best integrator from the computational accuracies and
efficiencies.

Now, we apply the optimized extended method N4P with the large
time step $\tau=0.1$ to explore the dynamics of the modified
H\'{e}non-Heiles system. As the magnitude of negative initial
value of $y$ decreases, there is a dynamical transition from
physical many-islands tori (Figs. 6 (a)-(d)), to single torus
(Fig. 6(e)) and to chaotic orbits (Fig. 6(f)). This seems to show
that the strength of chaos is enhanced with a decrease of the
magnitude of negative initial value of $y$. However, chaos is not
always enhanced, as can be seen from the dependence of FLI on the
initial value of $y$ in Fig. 7 that displays  the dynamical
transition from order to chaos. Chaos mainly occurs for the
initial values of $y$ in the vicinity of -2.25$\sim$-2.1, -1.6,
and  -1.2$\sim$-1. Here, the initial value $x=0$ is fixed, and
each value of FLI is obtained after integration time $t=3000$.
Numerical tests show that the threshold of FLIs between the
ordered and chaotic cases is 4. The FLIs larger than the threshold
determine the chaoticity of bounded orbits, but the FLIs less than
the threshold indicate the regularity of bounded orbits.

\subsection{Spring pendulum}

A spring pendulum in polar coordinates is described by the
Hamiltonian [1]
\begin{eqnarray}
  H=\frac{1}{2}\left(p_{r}^{2}+\frac{p_{\varphi}^{2}}{r^{2}}\right)-r \cos
  \varphi+(r-1)^{2}.
\end{eqnarray}
The Hamiltonian is divided into kinetic energy $K$ and potential
energy $V$ as follows:
\begin{eqnarray}
  K &=& \frac{1}{2}\left(p_{r}^{2} +\frac{p_{\varphi}^{2}}{r^{2}}\right), \\
  V &=& -r \cos \varphi+(r-1)^{2}.
\end{eqnarray}
$K$ is one of the expressions in Eq. (19) and is analytically
solved. In this case, the extended algorithms such as N4 can be
suitable for integrating the spring pendulum problem.

Taking the step size $\tau=0.1$ and energy $E=1/12$, we choose
initial conditions $r=1.15$, $p_{r}=0$ and $\varphi=0.05\pi$. The
initial value $p_{\varphi}>0$ is solved from the energy equation
$E=H$. Fig. 8 and Table 4 show that RK4 has the largest errors in
the  energy and position, and RKF8(9) exhibits the smallest energy
error. The errors in the energy have secular drifts for the
non-symplectic methods RK4 and RKF8(9), but do not have for the
symplectic methods M4, M4V, M4P, N4, N4O, N4V and N4P. The errors
in the  energy and position for N4V and N4P are about three orders
of magnitude smaller than those for M4. Table 5 lists CPU times of
the methods. The standard symplectic methods are slower than the
corresponding extended schemes in computational efficiencies. The
optimized methods  N4V and N4P need small additional cost compared
with the corresponding non-optimized method N4.

The tested orbit in Fig. 8 is a regular single closed torus on the
Poincar\'{e} section $\varphi=0$ and $p_{\varphi}>0$ in Fig 9(a).
Unlike in the modified H\'{e}non-Heiles system, M4 in the present
problem is the same as anyone of the methods M4V, M4P, N4, N4O,
N4V, N4P and RKF8(9) in the description of phase-space structures.
This is because the accuracies in the energy and position for M4
in Fig. 8 are higher than those in Fig. 1. In fact, the energy for
M4 is accurate to an order of $10^{-3}$, and the position for M4
is accurate to an order of 1 in Fig. 1 when the integration time
$t=10^4$. However, the energy for M4 is accurate to an order of
$10^{-5}$, and the position for M4 is accurate to an order of 0.1
in Fig. 8 when the integration time $t=10^4$. Because of this, M4
can truly describe the phase-space structures, as the methods M4V,
M4P, N4, N4O, N4V, N4P and RKF8(9) can. For given integrator and
time step, the numerical accuracy closely depends on model
Hamiltonians, and particularly depends on the periods of orbits in
the Hamiltonians. The larger the periods are, the better the
numerical accuracy is.

Considering the best performance in numerical accuracies and
computational efficiency, we use N4P to trace different
phase-space structures  when various initial values are given to
$\varphi$. For example, single-torus orbits exist for
$\varphi=0.19\pi$ in Fig. 9(b) and $\varphi=0.358\pi$ in Fig.
9(e). There are many islands for $\varphi= 0.35\pi$ in Fig. 9(d),
$\varphi=0.361\pi$ in Fig. 9(f), $\varphi=0.366\pi$ in Fig. 9(g),
and $\varphi=0.39\pi$ in Fig. 9(i). Chaos occurs for
$\varphi=0.2\pi$ in Fig. 9(c) and $0.376\pi$ in Fig. 9(h). No
universal rule for the dynamical transition from order to chaos
seems to be given to a variation of the initial value $\varphi$.
However, Fig. 10(a) for the description of the initial values of
$\varphi$ and their corresponding FLIs shows that chaos mainly
occurs for the initial values of $\varphi$ in the vicinity of
0.25, 1.75, and 2.25. This is because the spring pendulum suffers
from strong perturbations in the vicinity of $\varphi=\pi/4,
7\pi/4$. Scanning the initial values of $r$  and their
corresponding FLIs in  Fig. 10(b) displays that chaos mainly
occurs for the initial values of $r$ in the vicinity of 0.5, and
0.75$\sim$2.25.

Besides the method of FLIs, the 0-1 test for
chaos [33] can be applied to explore the transition from order to
chaos with the initial value $\varphi$ or $r$ varying. The 0-1
test chaos indicator is described here. Set
$\textbf{Z}=(r,\varphi,p_r, p_{\varphi})$ as a solution of the
system (39) at time $t$. $\psi(\textbf{Z})$ is a function of
$\textbf{Z}$; for example, $\psi(\textbf{Z})=r$ is a simple
choice. The authors of [33] defined two functions
\begin{eqnarray}
\theta(t) &=& ct+\int^{t}_{0}\psi(\textbf{Z}(s))ds,\\
q(t) &=& \int^{t}_{0}\psi(\textbf{Z}(s))\cos(\theta(s))ds,
\end{eqnarray}
where $c>0$ is an arbitrarily constant.  The mean-square
displacement of $q(t)$ is
\begin{eqnarray}
L(t)=\lim _{T \rightarrow \infty} \frac{1}{T}
\int_{0}^{T}[q(t+s)-q(s)]^{2} \mathrm{~d}s.
\end{eqnarray}
Finally, the asymptotic growth rate of the mean-square
displacement is expressed as
\begin{eqnarray}
\Lambda=\lim _{t \rightarrow \infty} \frac{\ln L(t)}{\ln t}.
\end{eqnarray}
Based on ergodic theory, $\Lambda=0$  signifies regular dynamics,
but $\Lambda=1$ signifies chaotic dynamics. Taking $c=1.8$,
$T=1000000$, $t=1000$ and time step $\tau=0.1$, we calculate the
values of $\Lambda$ with respect to the orbits in Fig. 10 and
obtain the dependence of $\Lambda$ on the initial values $\varphi$
in Fig. 11. The values of $\Lambda$ in Fig. 11(a) are in the
vicinity of 1 when the initial values of $\varphi$ are in the
vicinity of 0.25, 1.75 and 2.25, and are in the vicinity of 0 when
the initial values of $\varphi$ are in the vicinity of 0.1, 0.5,
1.5, 2.0 and 2.5. The values of $\Lambda$ in Fig. 11(b) are in the
vicinity of 1 when the initial values of $r$ are in the vicinity
of 0.5, and 0.75$\sim$2.25, and are in the vicinity of 0 when the
initial values of $r$ are in the vicinity of 0.7 and 2.5. That is
to say, the dynamical properties described by the 0-1 indicator
$\Lambda$ in Fig. 11 are consistent with those given by the FLIs
in Fig. 10. However, because $T$ is large enough, computations of
$L(t)$ are relatively expensive. In fact, CPU time is 284.41
seconds for each of the $\Lambda$ values in Fig. 11, and 0.36
seconds for each of the FLIs in Fig. 10. Thus, the FLIs are
quicker to distinguish between the ordered and chaotic two cases
than the 0-1 test indicator.

\section{Conclusions and discussions}

Many force-gradient explicit symplectic integration algorithms
with the force-gradient operator $C$ for the Hamiltonian (2) with
the kinetic energy (1) have been in Refs. [21-23]. However, these
algorithms become useless for the Hamiltonian (18) with the
integrable kinetic energy (19) if the force-gradient operator $C$
is not altered. We find that the existing integrators are still
available for the Hamiltonian (18) when the force-gradient
operator $C$ gives place to a new operator $D$. This new operator
is not a force-gradient operator but is similar to the
momentum-version operator associated to the potential $V$. The
extended algorithms are no longer solvers of the original
Hamiltonian but are solvers of slightly modified Hamiltonians.
They are explicit symplectic integrators with symmetry or
time-reversibility.

Numerical tests show that the standard symplectic integrators
without the operator $D$ are generally inferior to the
corresponding extended methods with the operator $D$ in
computational accuracies and efficiencies. For example, the
fourth-order Forest-Ruth symplectic scheme cannot provide reliable
results to the description of regular and chaotic dynamical
features of the modified H\'{e}non-Heiles system for the use of
some appropriately large time steps, but the corresponding
extended ones can. The optimized methods have better accuracies
than the corresponding non-optimized methods. Among the tested
symplectic methods, the two extended optimized seven-stage
fourth-order methods of Omelyan, Mryglod and Folk (N4V and N4P)
exhibit the best numerical performance, and their accuracies are
about three orders of magnitude better than the accuracies of the
Forest-Ruth symplectic scheme. Finally, one of the two optimized
algorithms is used to study the orbital dynamical features of the
modified H\'{e}non-Heiles system and the spring pendulum.

The proposed extended algorithms are suitably
applicable to any Hamiltonian systems like the Hamiltonian system
(18). The kinetic energies like Eq. (19) can be found in some
references. For example, the Hamiltonian in Eq. (5.21) of Ref.
[34] is
\begin{equation}
H=\frac{1}{2}(p^2_1+p^2_2)+\frac{\varepsilon}{2}p^2_1\cos x_2+
\cos x_1-1.
\end{equation}
The kinetic energies given in Eq. (40) are universal for
two-dimensional Hamiltonian systems
$H=m(\dot{x}^2+\dot{y}^2)/2+V(x,y)$ in polar coordinates
$(r,\varphi)$, where potentials $V(x,y)$ are non-axisymmetric. In
addition, three-dimensional Hamiltonian problems
$H=m(\dot{x}^2+\dot{y}^2+\dot{z}^2)/2+V(x,y,z)$ in spherical
coordinates $(r,\theta,\varphi)$ are expressed as
\begin{equation}
H=\frac{1}{2m}(p^2_r+\frac{p^2_{\theta}}{r^2}+\frac{p^2_{\varphi}}{r^2\sin^2\theta})+V(r,\theta,\varphi),
\end{equation}
which resembles the Hamiltonian (18) with the integrable kinetic
energy (19). In other words, the Hamiltonian (18) is not
restricted to several special examples, but is often met in many
situations. In this sense, this extension has wide applications.
As an example, the kinetic energies in the Hamiltonians of the
spiral structure in self-consistent models of rotating galaxies
[35] and the spiral arms in galaxies [36] are Eq. (40). Thus, the
newly extended force-gradient explicit symplectic methods should
be suitable for the study of chaotic spiral arms in the models of
rotating galaxies. This problem will be considered in a future
work.

\section*{Acknowledgments}

The authors are very grateful to three referees for valuable
comments and useful suggestions. This research has been supported
by the National Natural Science Foundation of China [Grant Nos.
11973020 (C0035736), and 12133003], the Special Funding for
Guangxi Distinguished Professors (2017AD22006), and the National
Natural Science Foundation of Guangxi (Nos. 2018GXNSFGA281007 and
2019JJD110006).

\newpage

\begin{table*}[htbp] \centering \caption{Energy
errors for the algorithms in Figs. 1 (a) and (c). Because RKF8(9)
uses variable time steps, $\tau$=0.1 and 0.01 are not time steps
but are time intervals for outputting data. Note that 1.29 denotes
the error with an order of $10^{1.29}$, -9.69 means the error with
an order of $10^{-9.69}$, and so on. } \label{Tab1}
\begin{tabular}{cccccccccccc}
\hline Method  & RK4  & RKF8(9) & M4 & M4P & M4V & N4 & N4O & N4P & N4V\\
\hline $\tau=0.1$  & 1.29 & -9.69 & -2.73 &-4.08 & -4.13 & -3.96 & -4.40 & -5.75 & -5.66 \\
\hline $\tau=0.01$ & -3.63 & -11.67 & -6.75 & -8.09 & -8.14 & -7.97 & -8.40 & -9.72 & -9.67 \\
\hline
\end{tabular}
\end{table*}

\begin{table*}[htbp] \centering \caption{ Position
errors for the algorithms in Figs. 1 (b) and (d). } \label{Tab2}
\begin{tabular}{cccccccccccc}
\hline Method   & RK4  & M4 & M4P & M4V & N4 & N4O & N4P & N4V\\
\hline $\tau=0.1$  & 0.47 & 0.006 & -0.56  & -0.63 & -0.49 & -0.87 & -2.06 & -2.03 \\
\hline $\tau=0.01$ & -2.38 & -2.32 & -4.07 & -4.50 & -3.96 & -4.72 & -5.858 & -5.856 \\
\hline
\end{tabular}
\end{table*}

\begin{table*}[htbp] \centering
\caption{CPU times (unit: second) of each algorithm in Figs. 1 (a)
and (c). } \label{Tab3}
\begin{tabular}{cccccccccccc}
\hline Method   & RK4  & RKF8(9) & M4 & M4P & M4V & N4 & N4O & N4P & N4V\\
\hline $\tau=0.1$  & 0.33 & 24 & 0.51 & 0.66  & 0.64 & 0.44 & 0.44 & 0.59 & 0.59 \\
\hline $\tau=0.01$ & 2.56 & 83 & 4.53 & 5.89 & 5.70 & 3.63 & 3.64 & 5.22 & 5.23 \\
\hline
\end{tabular}
\end{table*}

\begin{table*}[htbp] \centering \caption{ Energy
errors $|\Delta H|$ and position errors $|\Delta r|$ for the
algorithms in Fig. 8.  } \label{Tab4}
\begin{tabular}{cccccccccccc}
\hline Method   & RK4  & RKF8(9) & M4 & M4P & M4V & N4 & N4O & N4P & N4V\\
\hline $|\Delta H|$  & 0.04 & -10.53 & -4.47 & -5.73  & -5.65 &-5.73 & -5.74 & -7.65 & -7.47 \\
\hline $|\Delta r|$  & 0.13 &  & -0.67 &-2.99 & -2.74 & -3.06 & -3.45 & -4.34 & -4.24 \\
\hline
\end{tabular}
\end{table*}

\begin{table*}[htbp] \centering \caption{ CPU
times (unit: second) of each algorithm in Fig. 8(a).  }
\label{Tab5}
\begin{tabular}{cccccccccccc}
\hline Method   & RK4  & RKF8(9) & M4   & M4P   & M4V  & N4   & N4O  & N4P  & N4V\\
\hline Time     & 0.73 & 28      & 3.70 & 4.56  & 3.86 & 2.78 & 2.98 & 3.89 & 3.30 \\
\hline
\end{tabular}
\end{table*}

\newpage

\begin{figure*}
\center{
\includegraphics[scale=0.35]{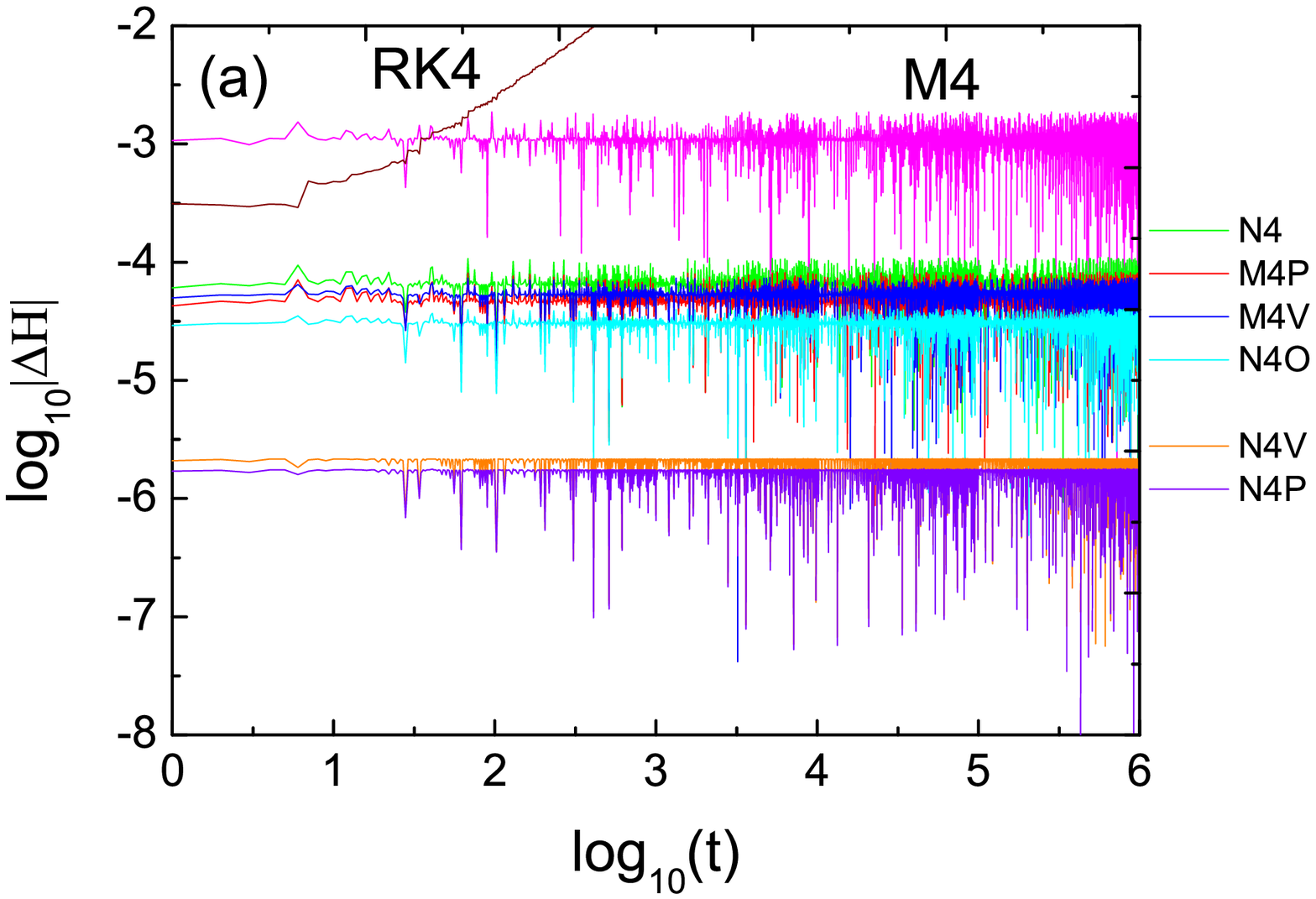}
\includegraphics[scale=0.35]{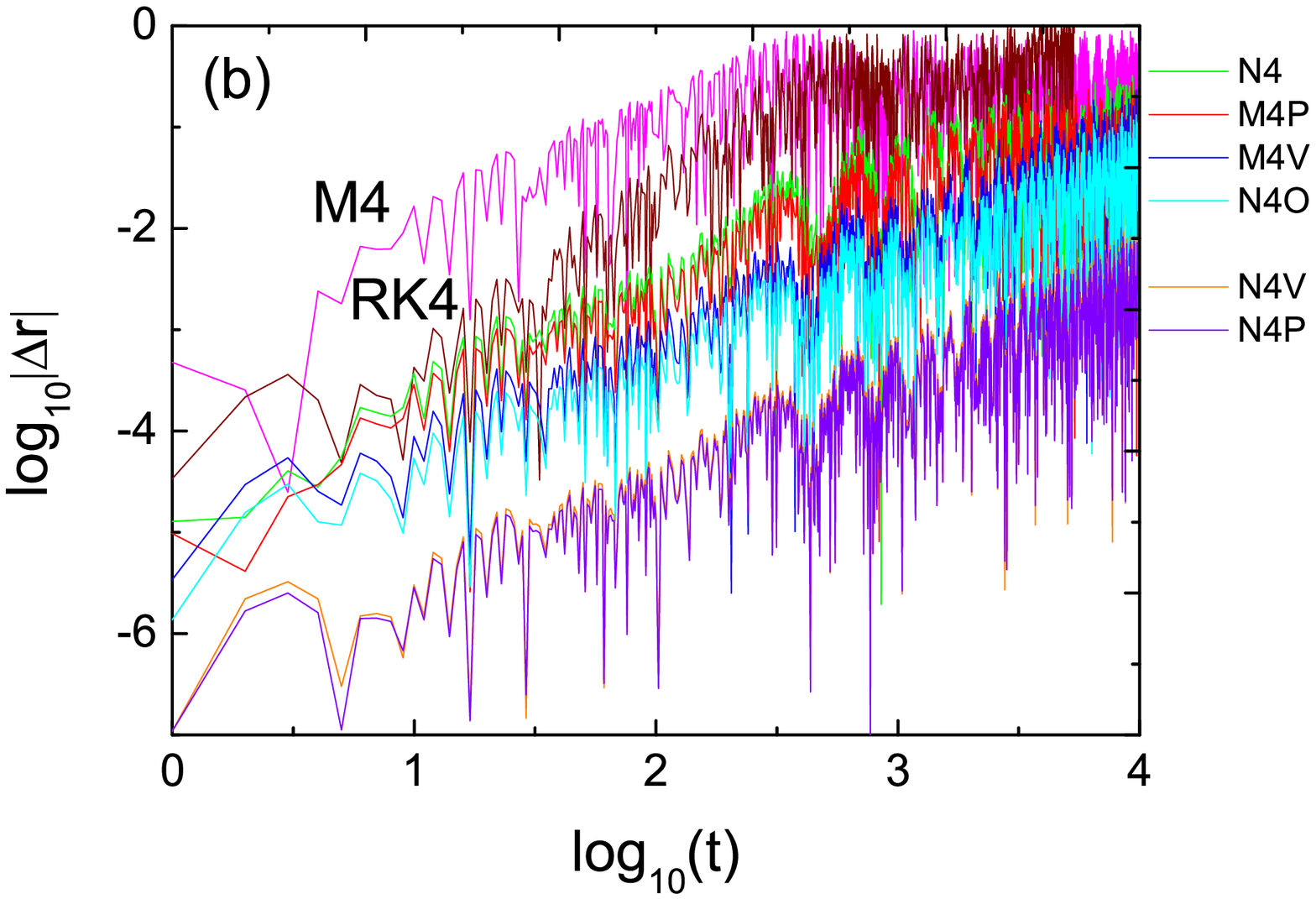}
\includegraphics[scale=0.35]{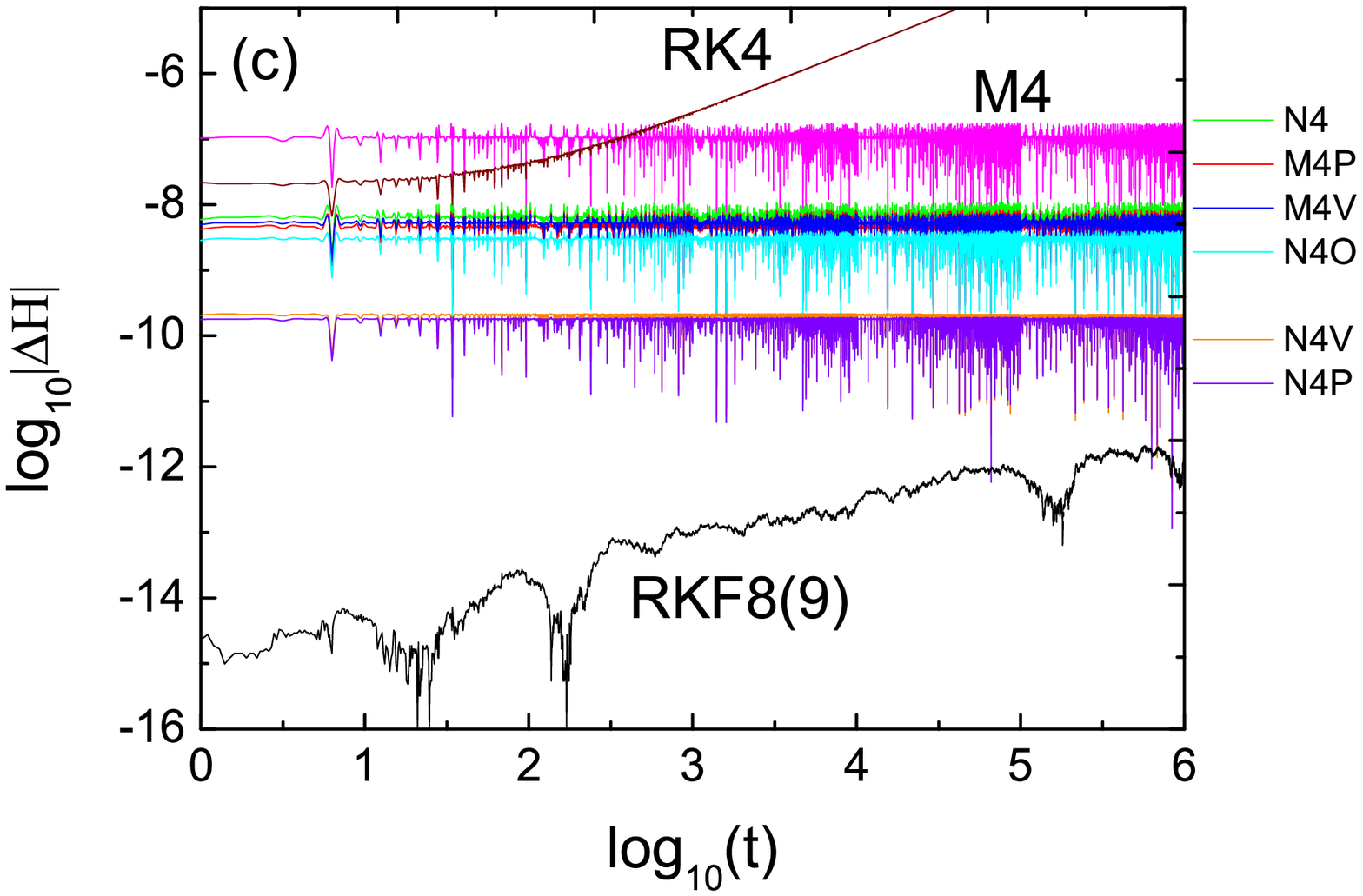}
\includegraphics[scale=0.35]{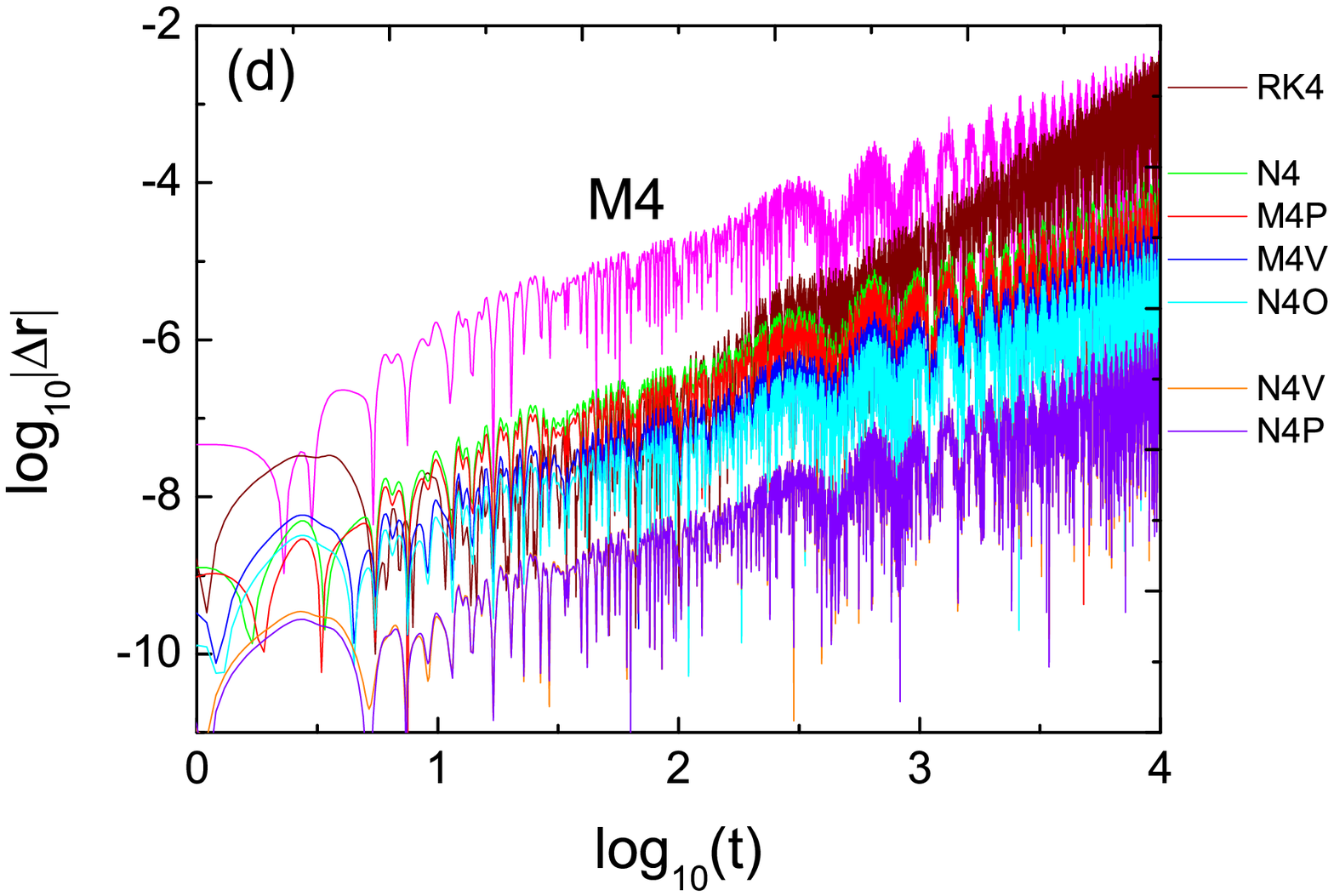}
\caption{(a) Energy errors for several algorithms independently
solving  the modified H\'{e}non-Heiles system. The time step is
$\tau=0.1$, and energy is $E=1/120$. Initial conditions are $x=0$,
$y=-2.02$ and $p_{y}=0$. The positive initial value of $p_x$ is
given by the energy relation $E=H$.  (b) Position errors between
the solutions of RKF8(9) and the other methods. (c) and (d) are
the same as (a) and (b) but a smaller time step $\tau=0.01$ is
used. These errors are clearly listed in Tables 1 and 2.  } }
\label{f1}
\end{figure*}

\begin{figure*}
\center{
\includegraphics[scale=0.5]{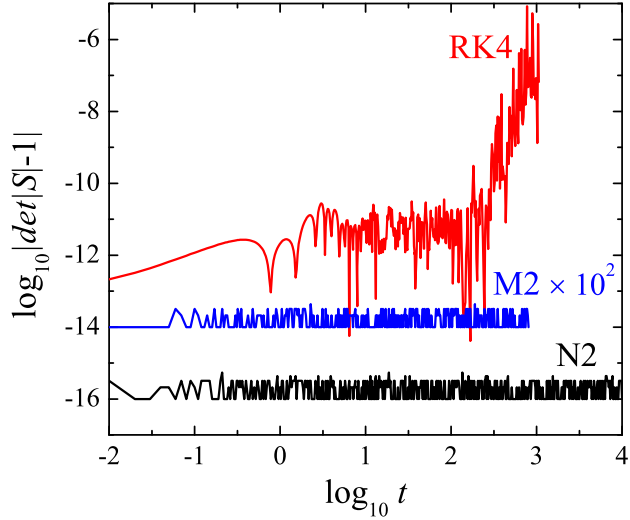}
\caption{Difference between 1 and the determinant
$det|\textbf{S}|$ of the matrix $\textbf{S}$ for each of the three
methods RK4, M2 and N2. The time step is $\tau=0.01$ and energy is
$E=1/120$. The initial conditions are $x=0$, $y=-0.988$ and
$p_{y}=0$. $M2\times 100$ means that the plotted values are 100
times larger than the practical values.  } } \label{f2}
\end{figure*}

\begin{figure*}
\center{
\includegraphics[scale=0.25]{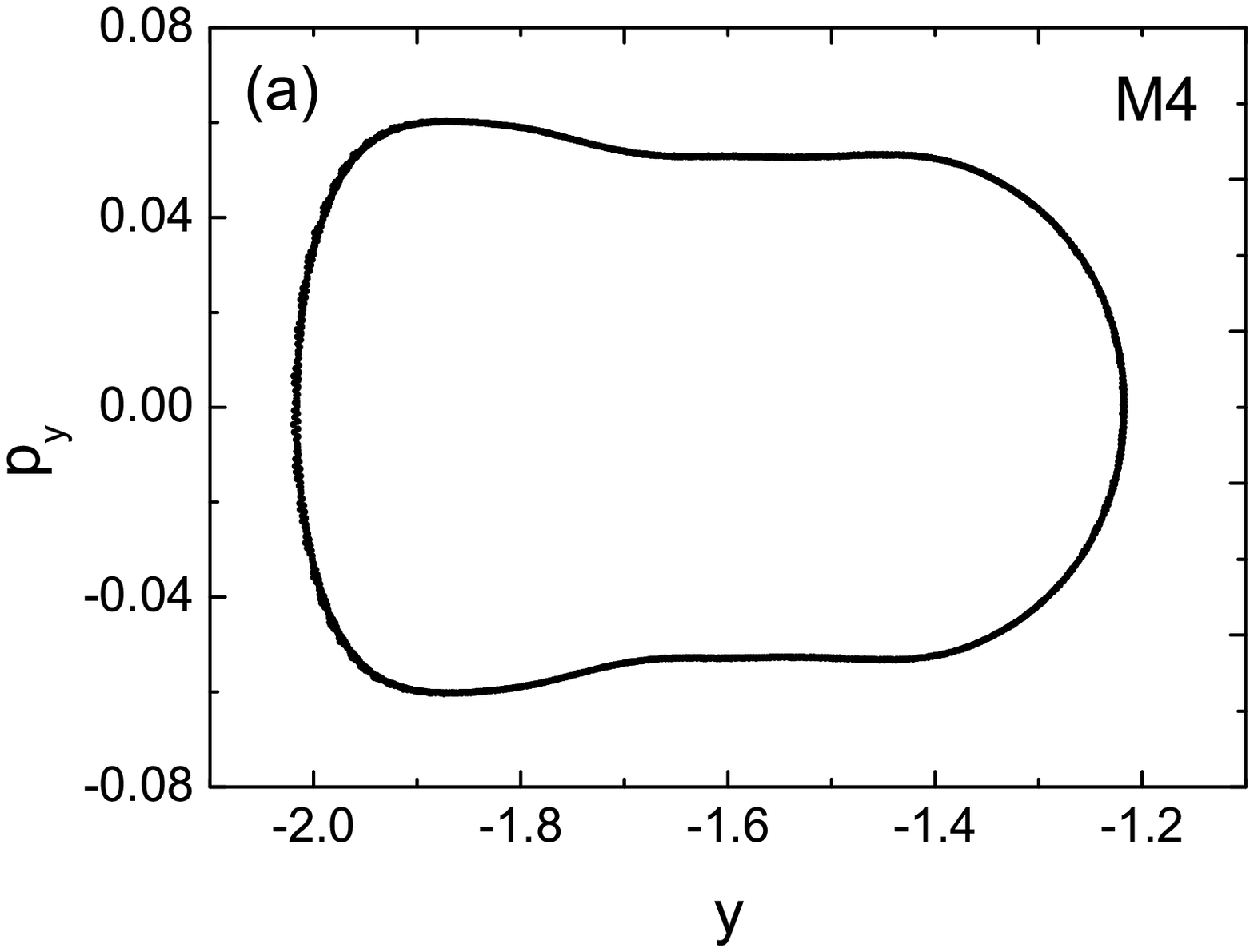}
\includegraphics[scale=0.25]{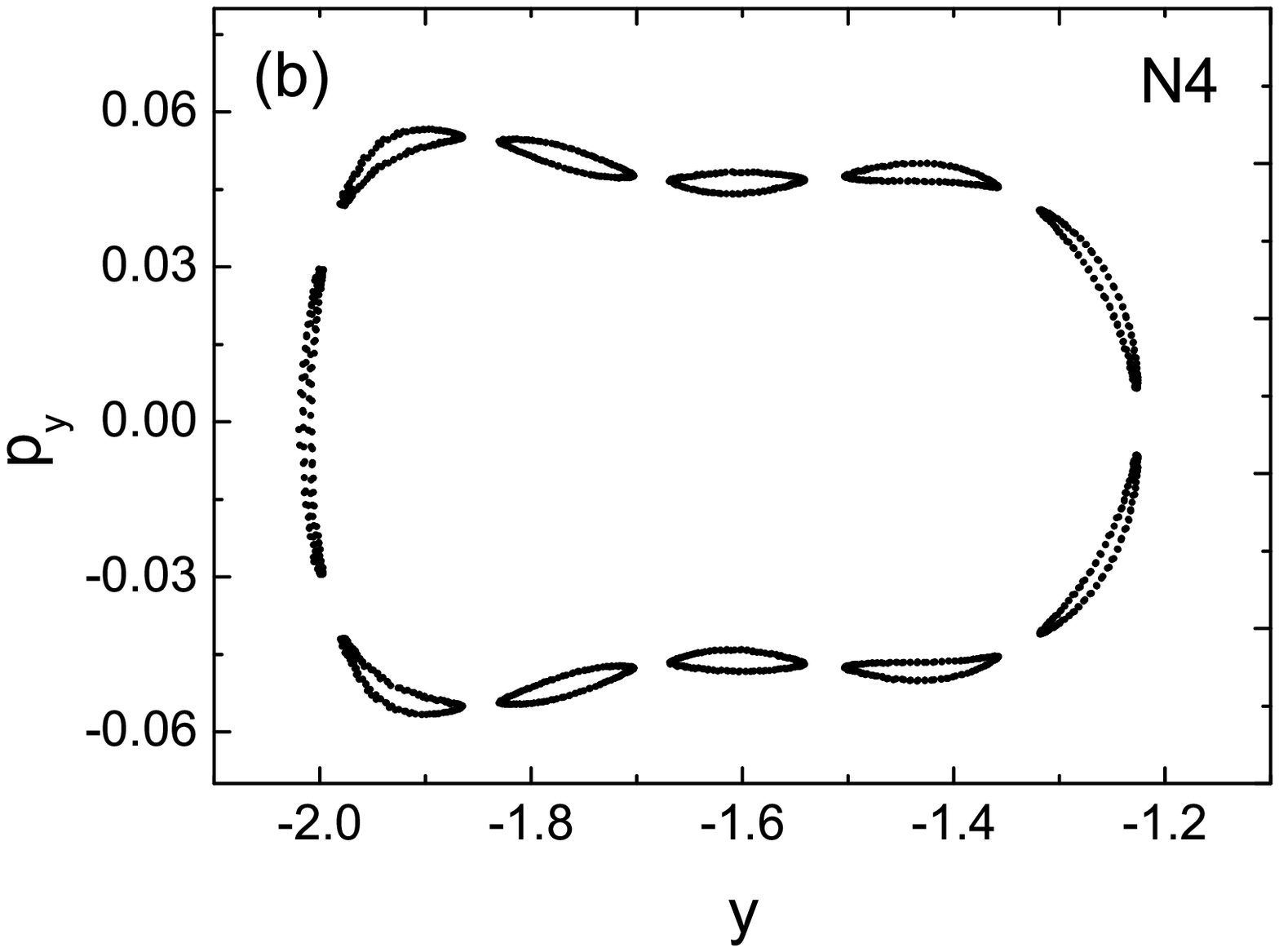}
\includegraphics[scale=0.25]{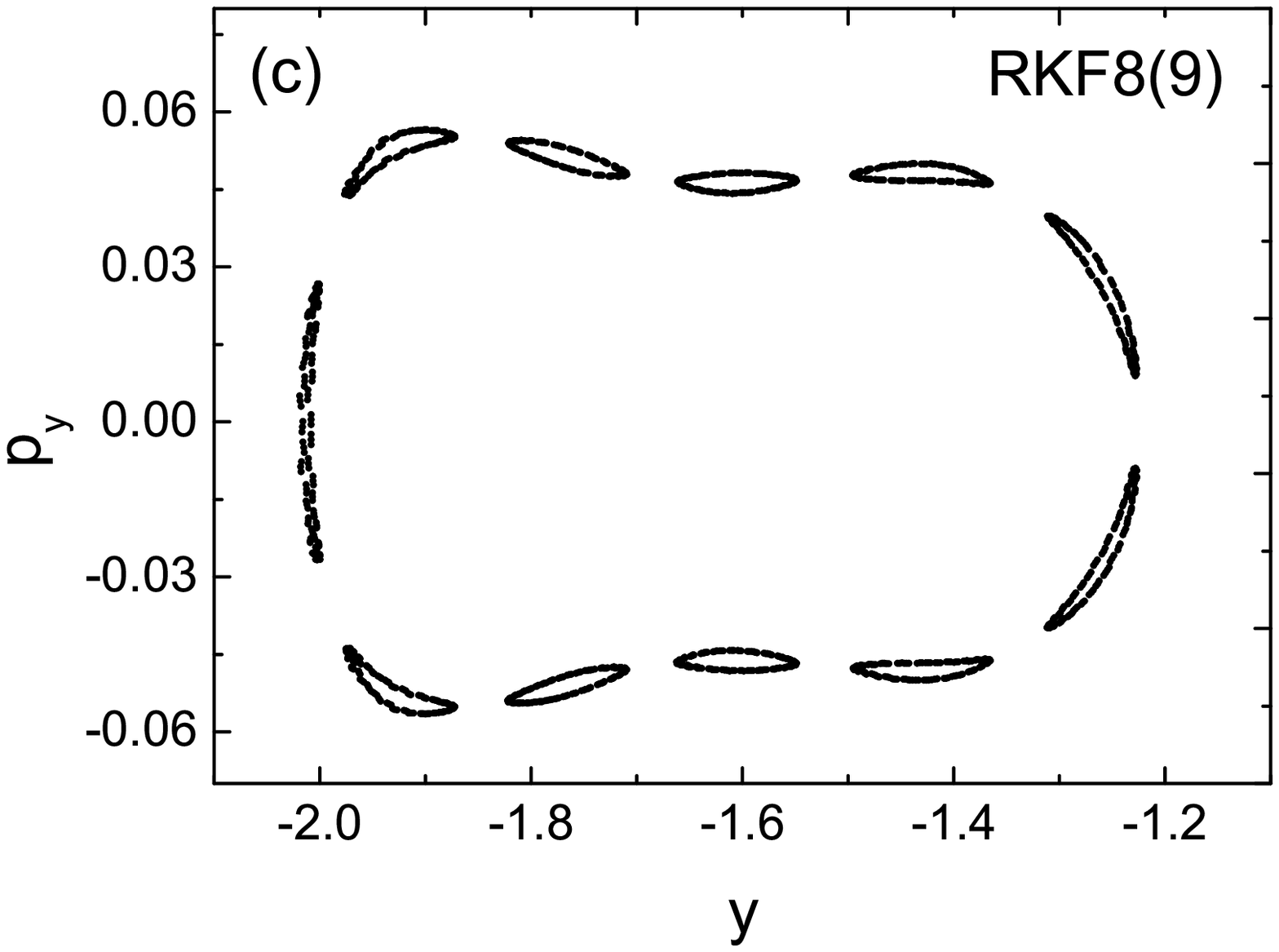}
\caption{Poincar\'{e} section $x=0$ and $p_{x}>0$. The integrators
M4 and N4 with the larger time step $\tau=0.1$ are respectively
used to act on the orbit tested in Fig. 1. M4 gives a single torus
to the orbit, but N4 like RKF8(9) exhibits many islands. The
result is reliable for N4, but is not for M4 because the larger
time step $\tau=0.1$ causes M4 to yield poor accuracy to the
numerical solutions. }} \label{f3}
\end{figure*}

\begin{figure*}
\center{
\includegraphics[scale=0.35]{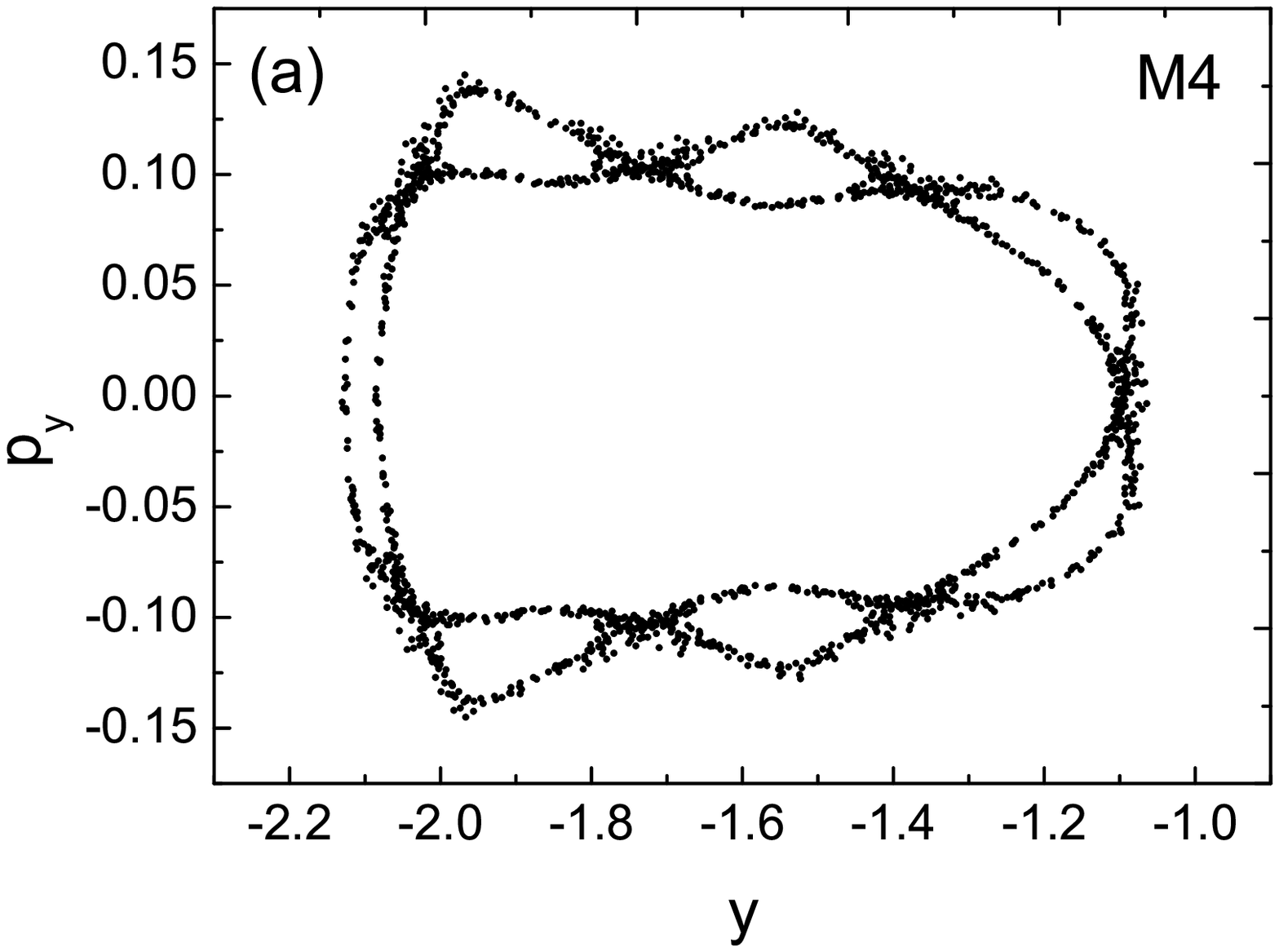}
\includegraphics[scale=0.35]{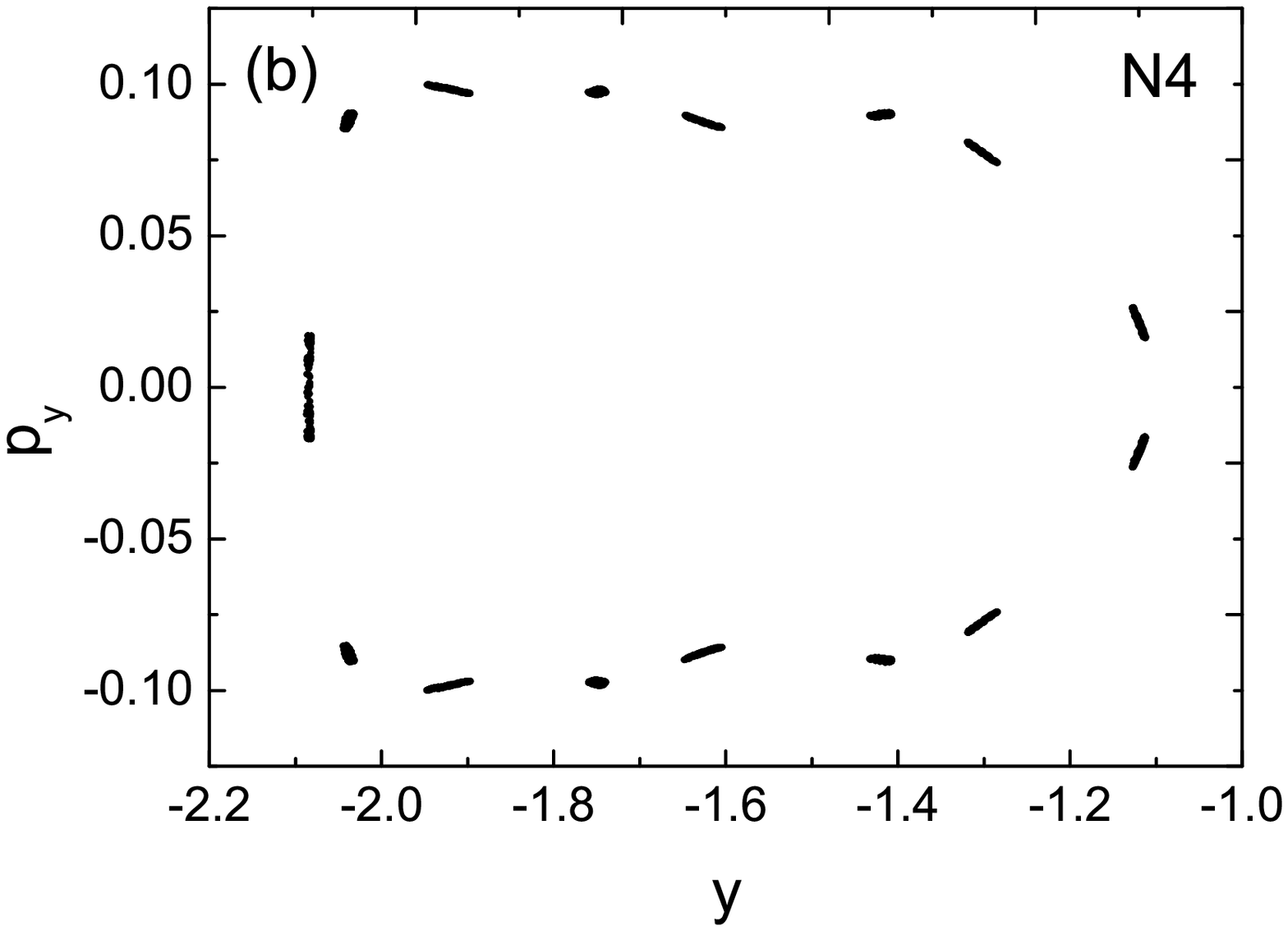}
\includegraphics[scale=0.35]{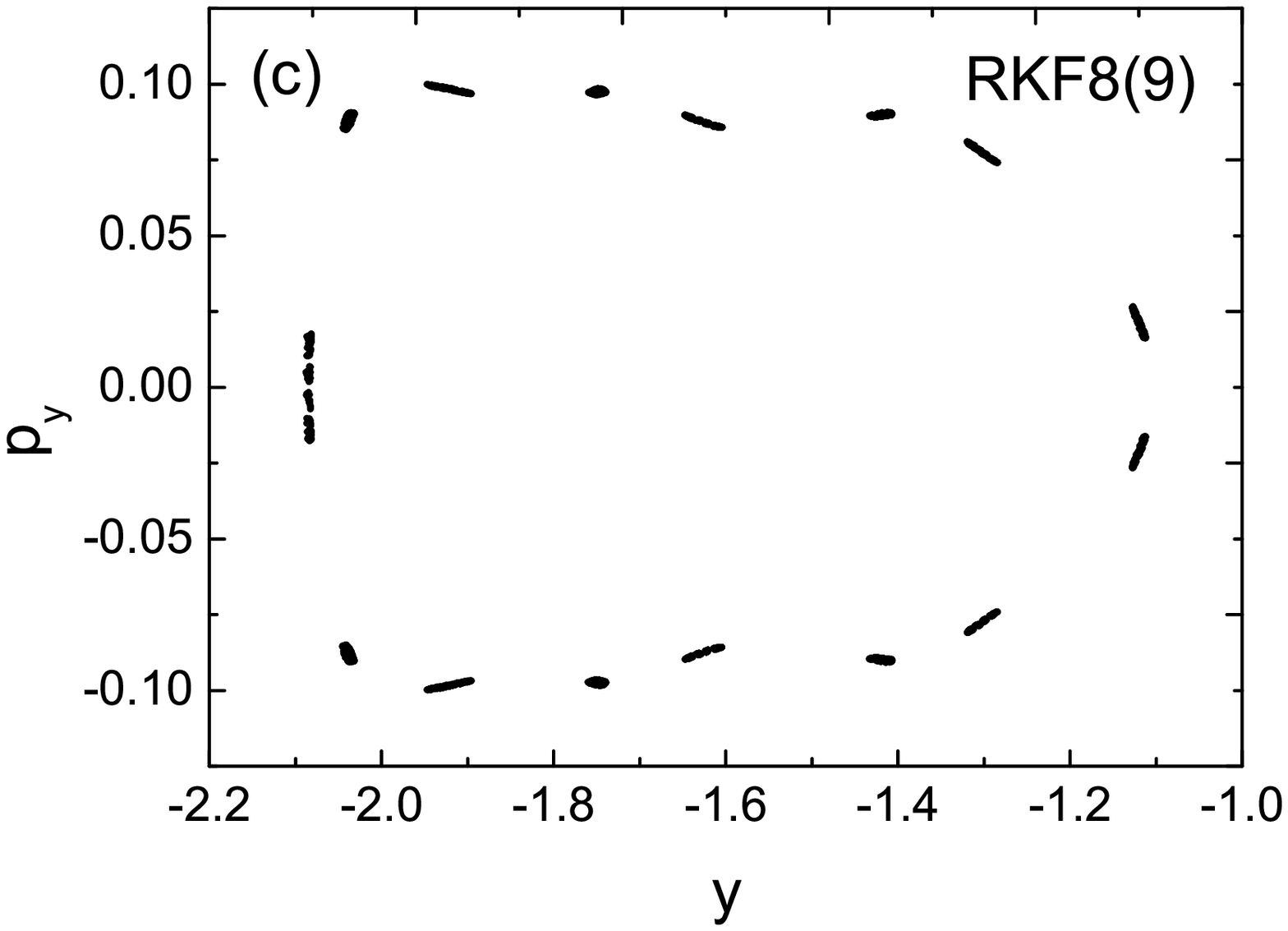}
\includegraphics[scale=0.35]{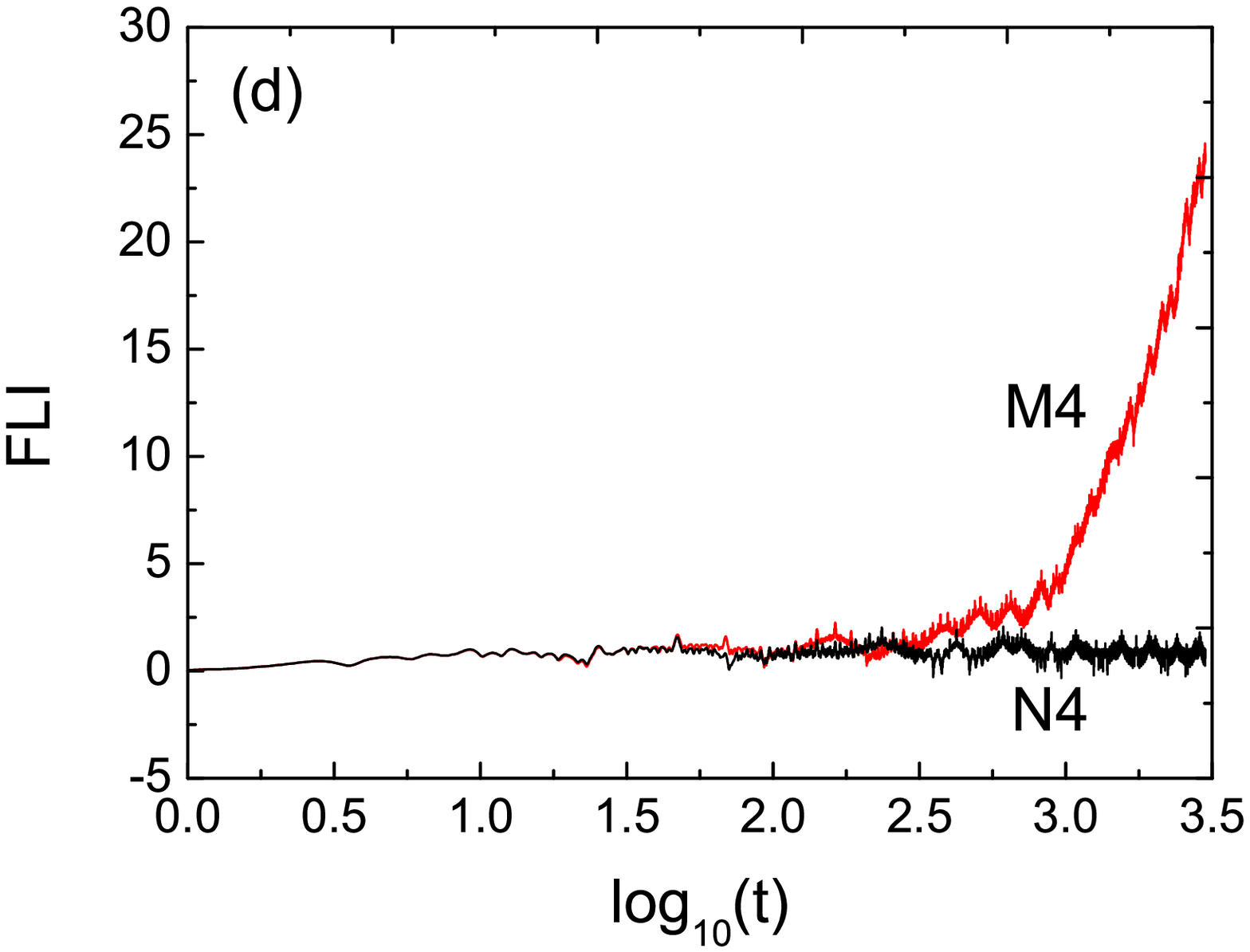}
\caption{(a)-(c): Same as Fig. 3 but for different initial values
$y=-1.108$ and $ p_{x}$. The orbit for M4 in (a) seems to be a
chaotic many-islands torus, but the orbits for N4 and RKF8(9) in
(b) and (c) are regular many-islands tori only. (d) Fast Lyapunov
indicators (FLIs) for the two methods M4 and N4 integrating the
orbit. The FLIs given by RKF8(9) are almost the same as those
given by N4. When the integration ends, the FLI is smaller than
2.5 for N4, but reaches 25 for M4. The FLI for N4 indicates the
regularity, while the FLI for M4 shows the chaoticity. The results
are correct for N4, whereas chaos for M4 is spurious in panels (a)
and (d) due to the large time step $\tau=0.1$ leading to M4 with
poor accuracy to the numerical solutions.  }} \label{f4}
\end{figure*}

\begin{figure*}
\center{
\includegraphics[scale=0.35]{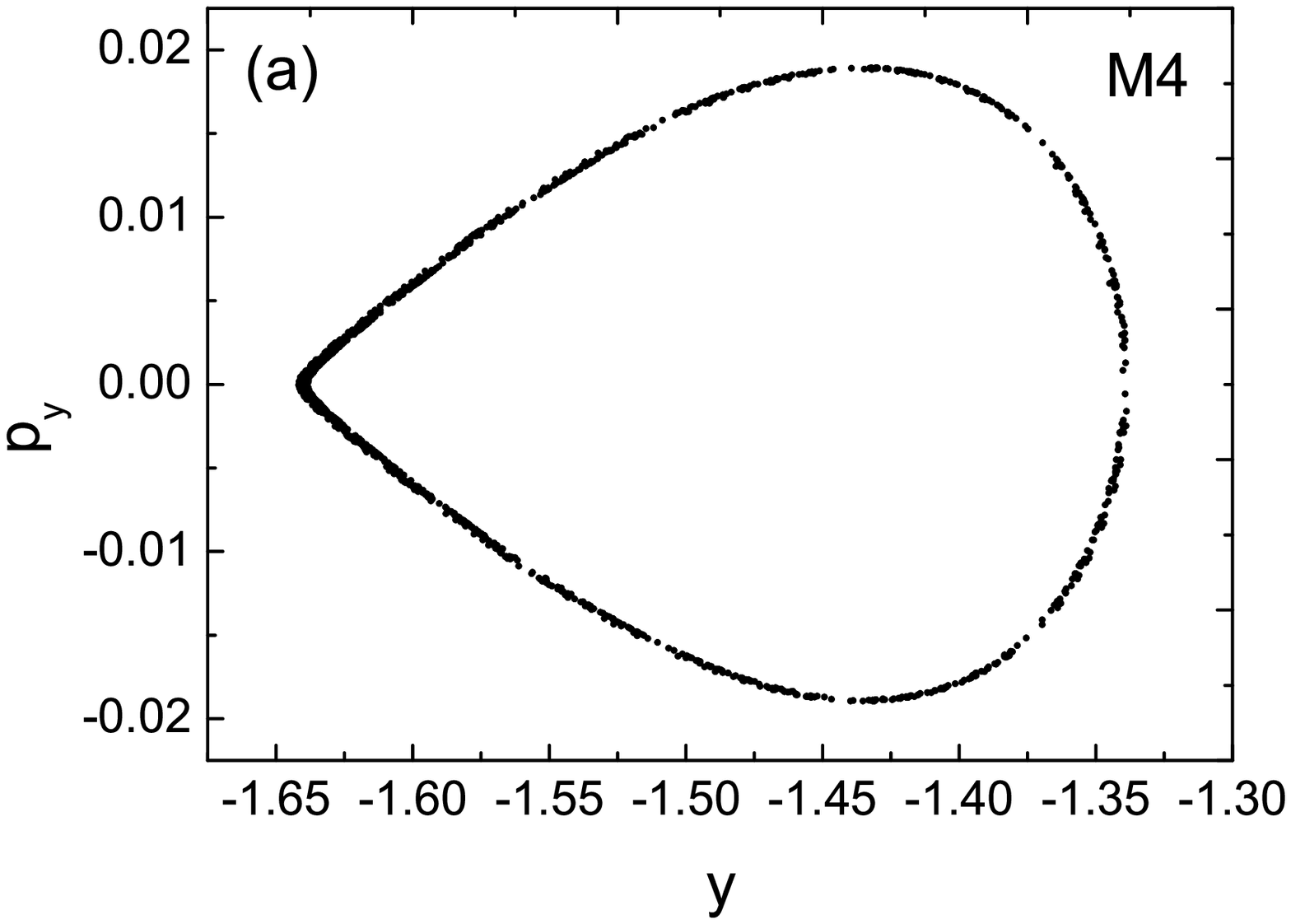}
\includegraphics[scale=0.35]{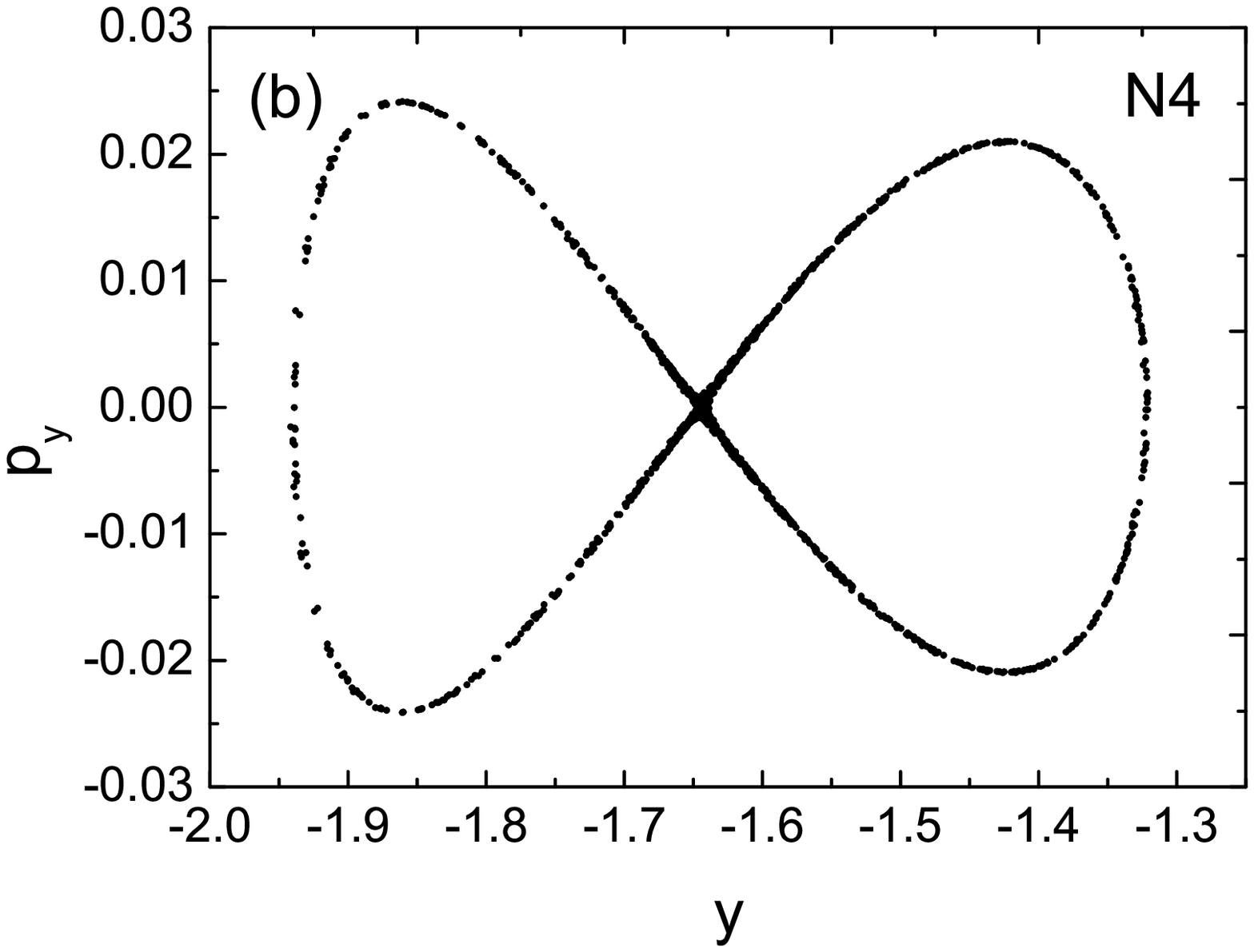}
\includegraphics[scale=0.35]{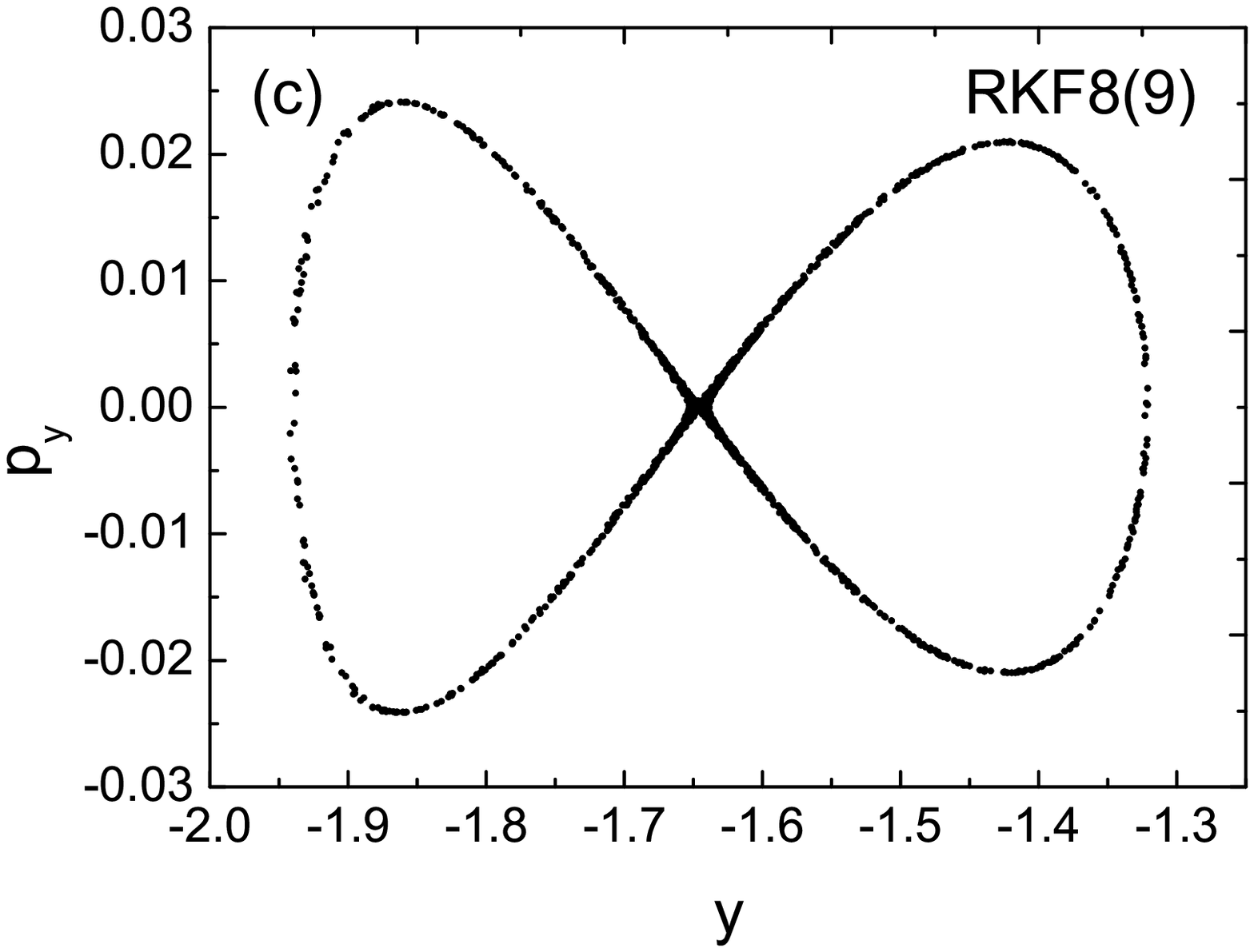}
\includegraphics[scale=0.35]{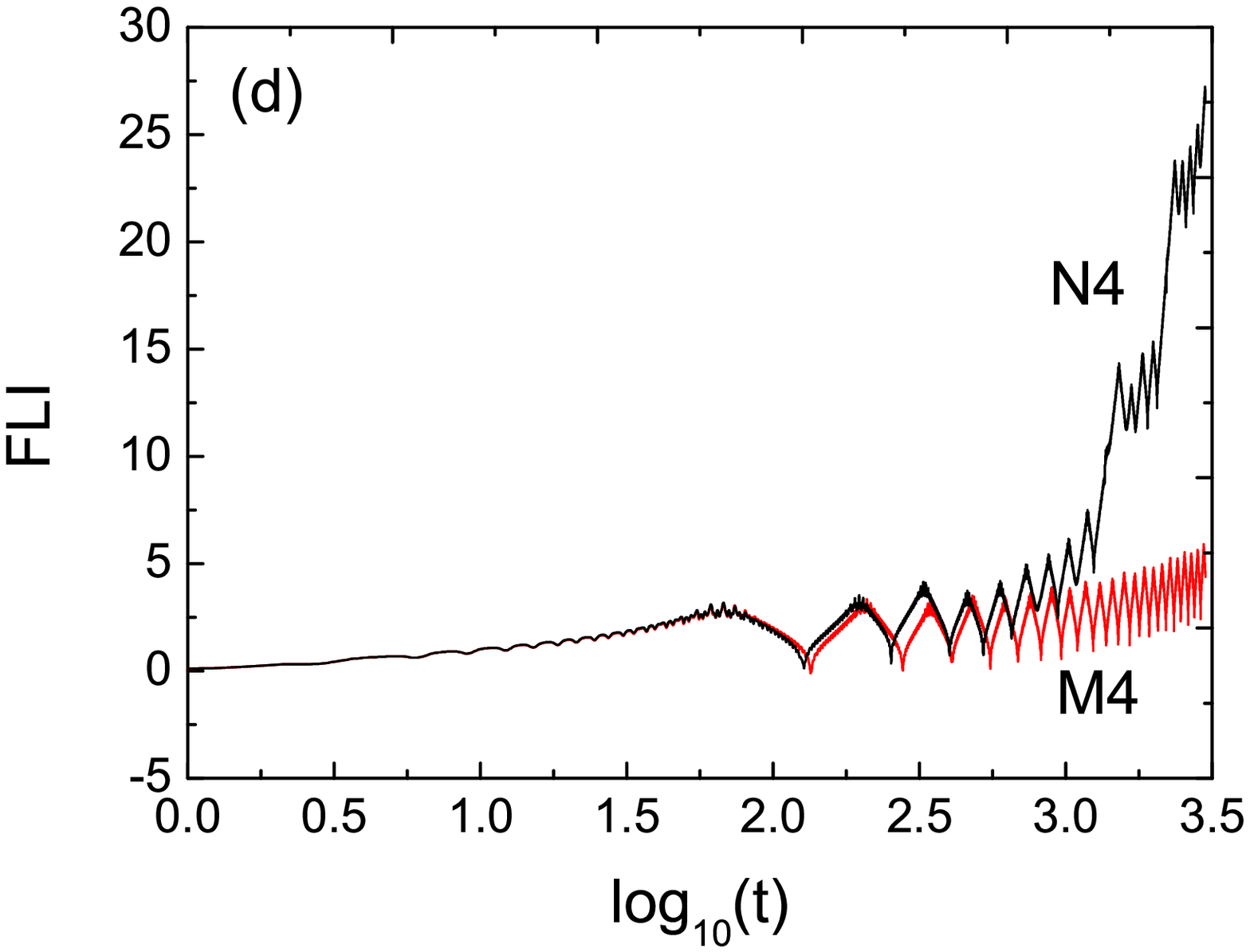}
\caption{Same as Fig. 4 but for different initial values
$y=-1.654$ and $ p_{x}$. (a): Regular single torus. (b) and (c):
Figure-eight orbits seem to be regular, but are chaotic, as shown
by the FLIs in (d). The regularity for M4 in panels (a) and (d) is
spurious, but the chaoticity for N4 is physical because N4
provides more accurate numerical solutions than M4.  }} \label{f5}
\end{figure*}

\begin{figure*}
\center{
\includegraphics[scale=0.26]{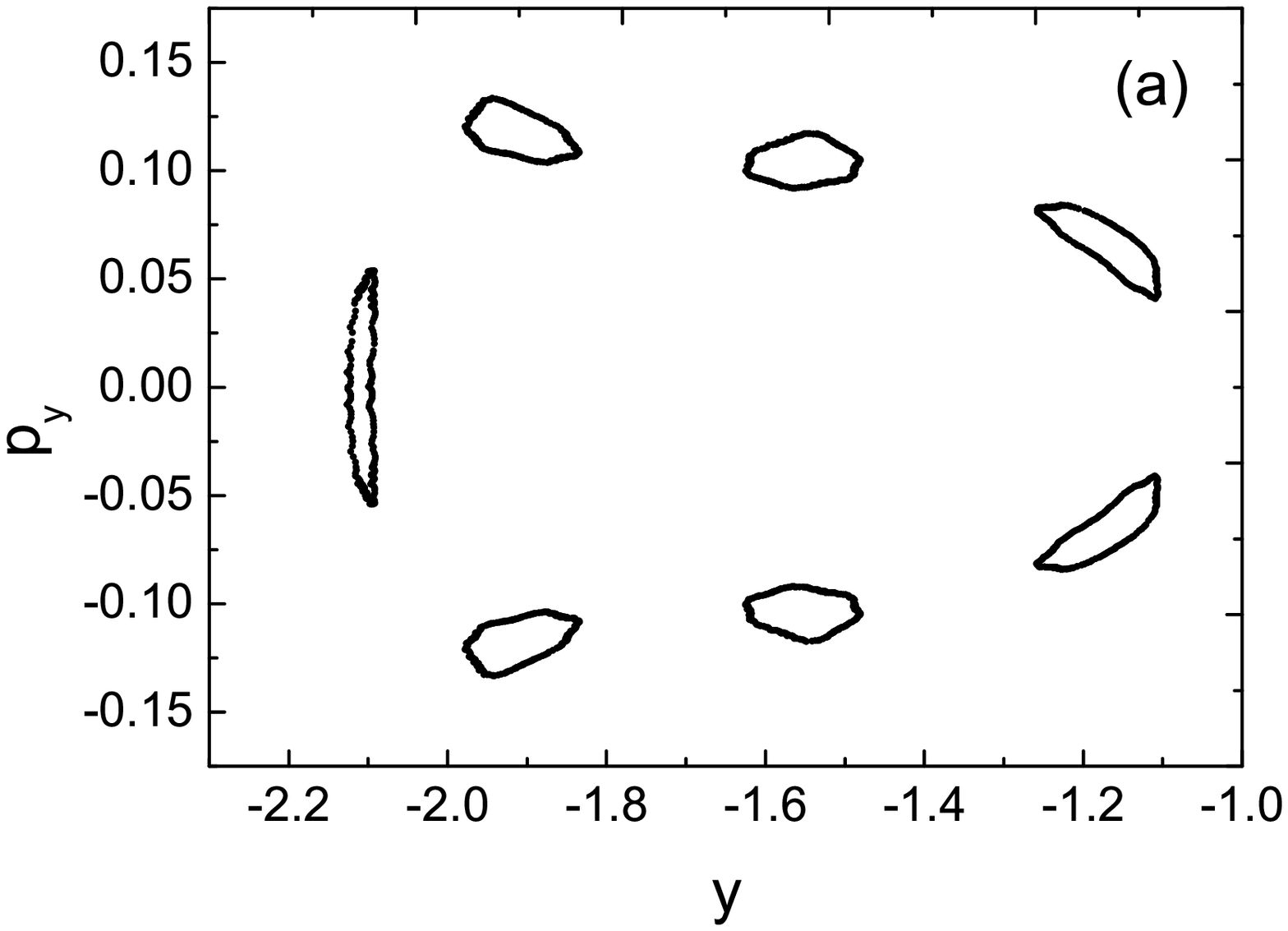}
\includegraphics[scale=0.26]{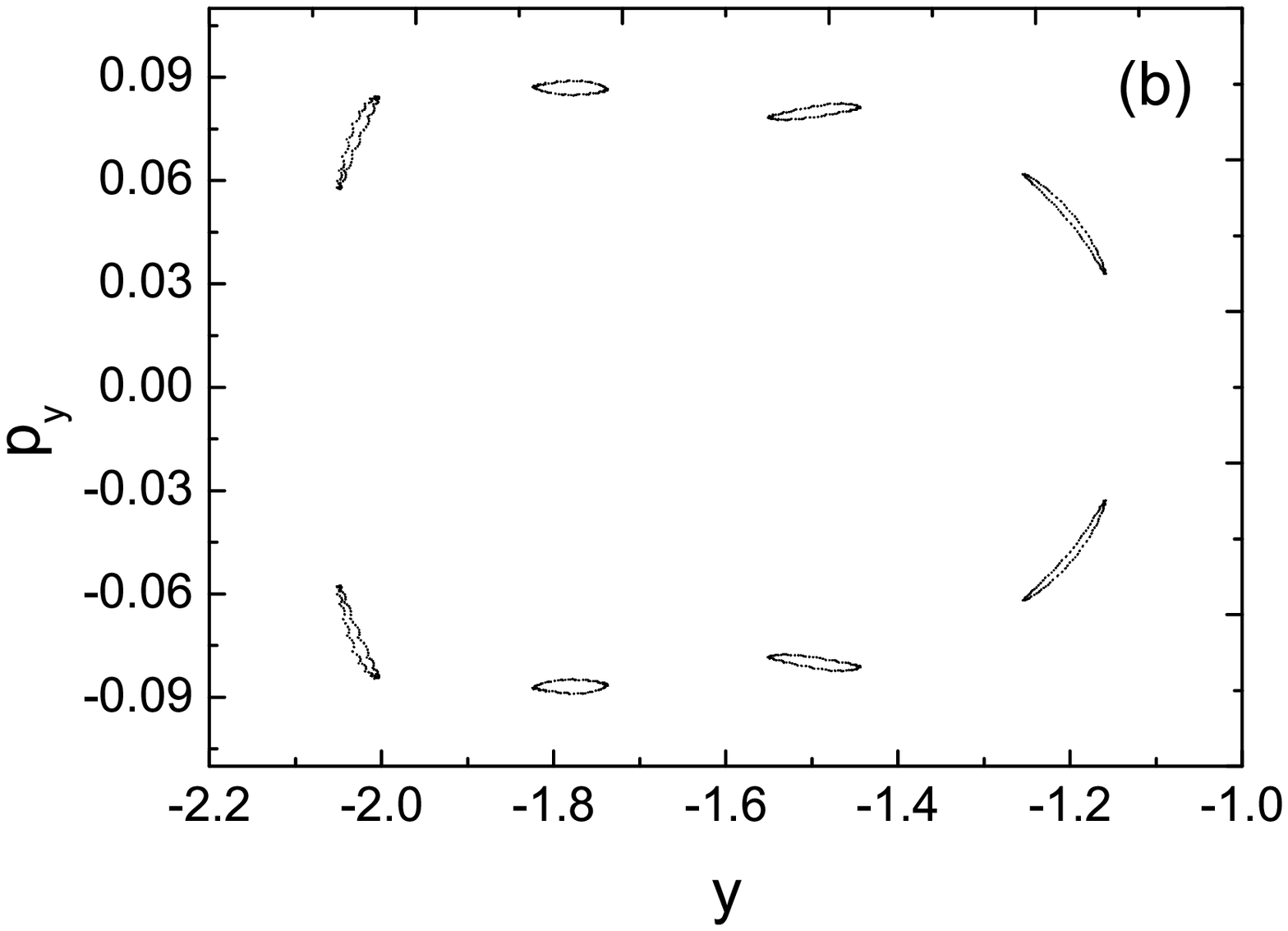}
\includegraphics[scale=0.26]{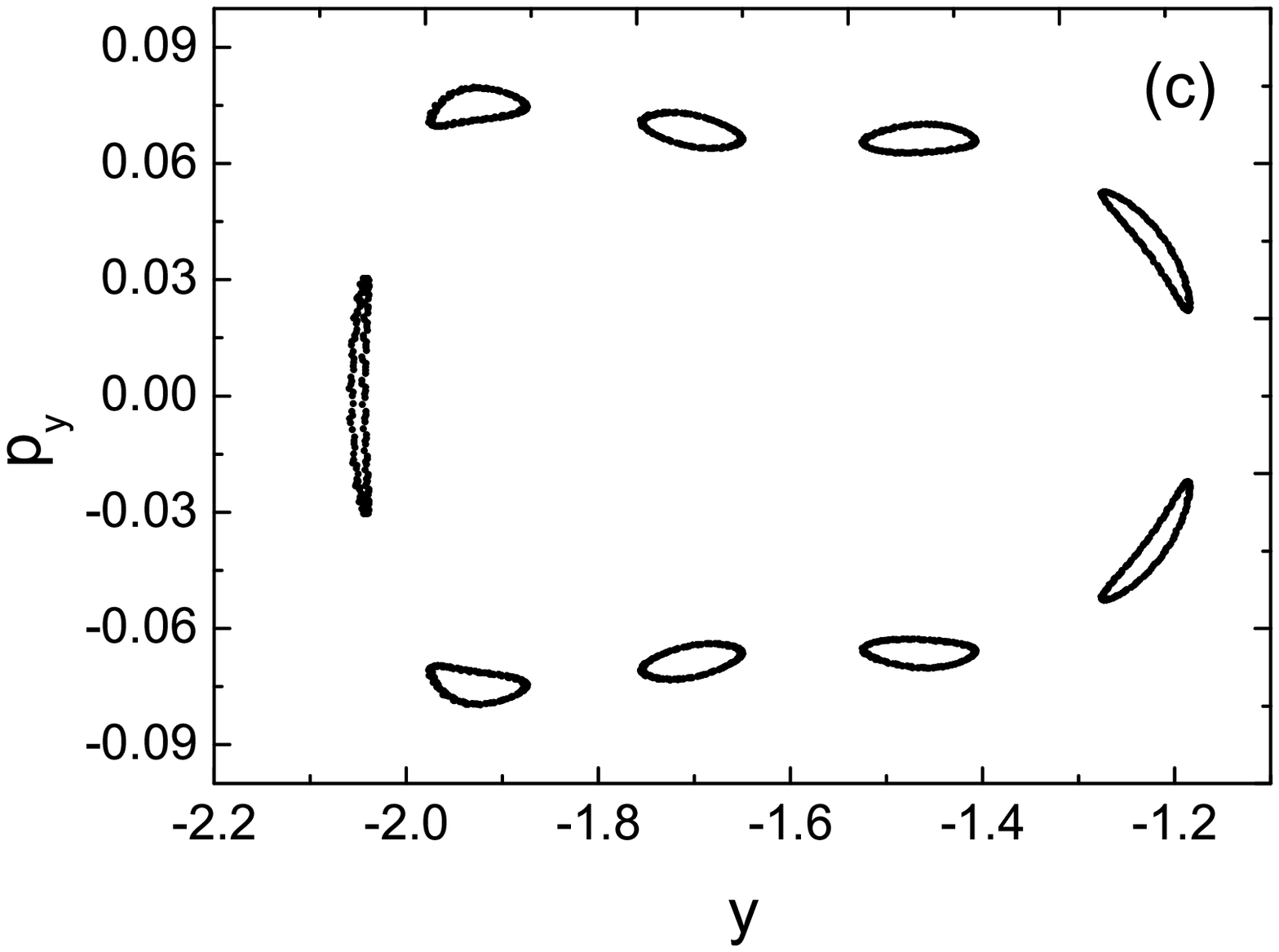}
\includegraphics[scale=0.26]{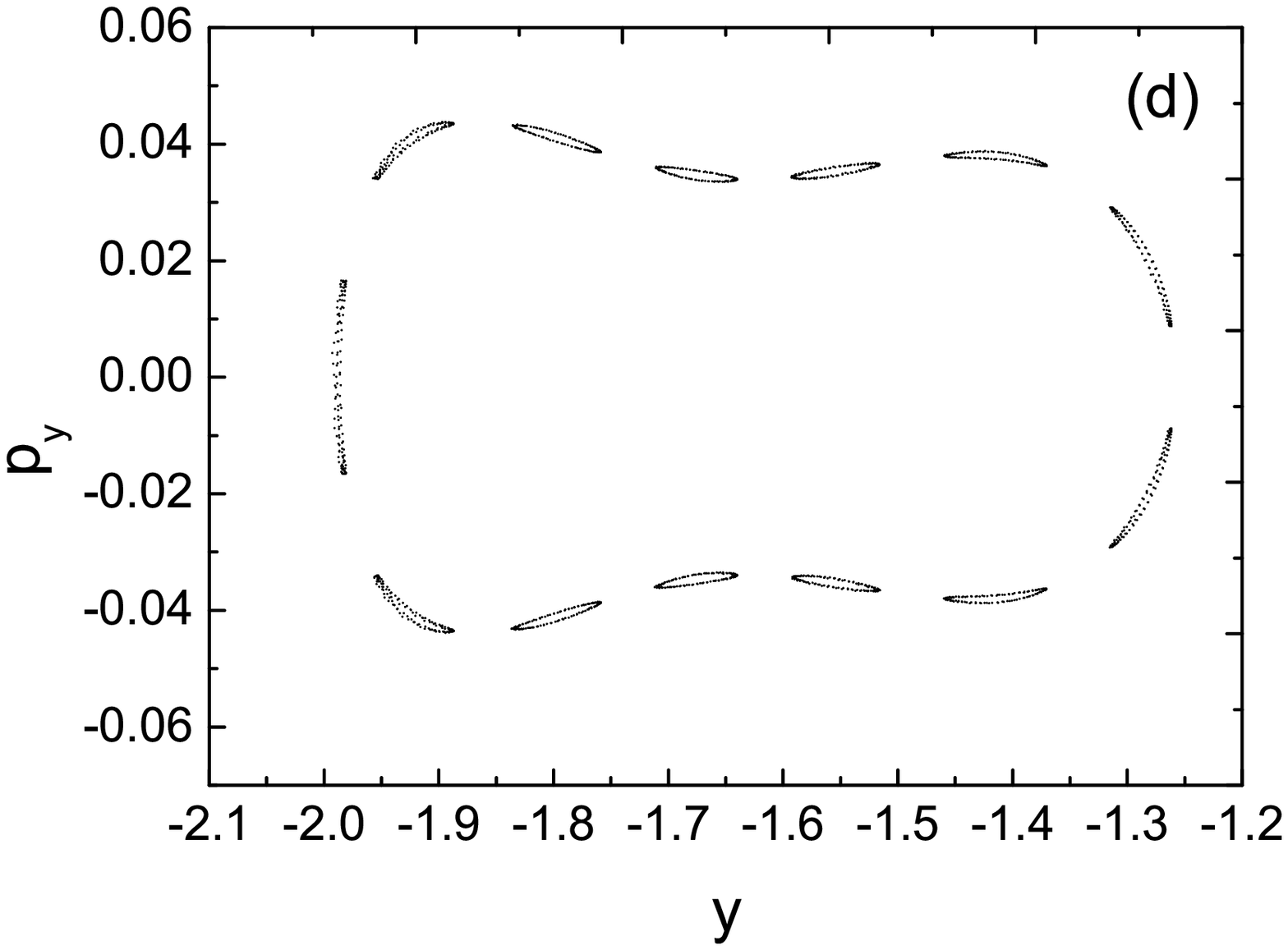}
\includegraphics[scale=0.26]{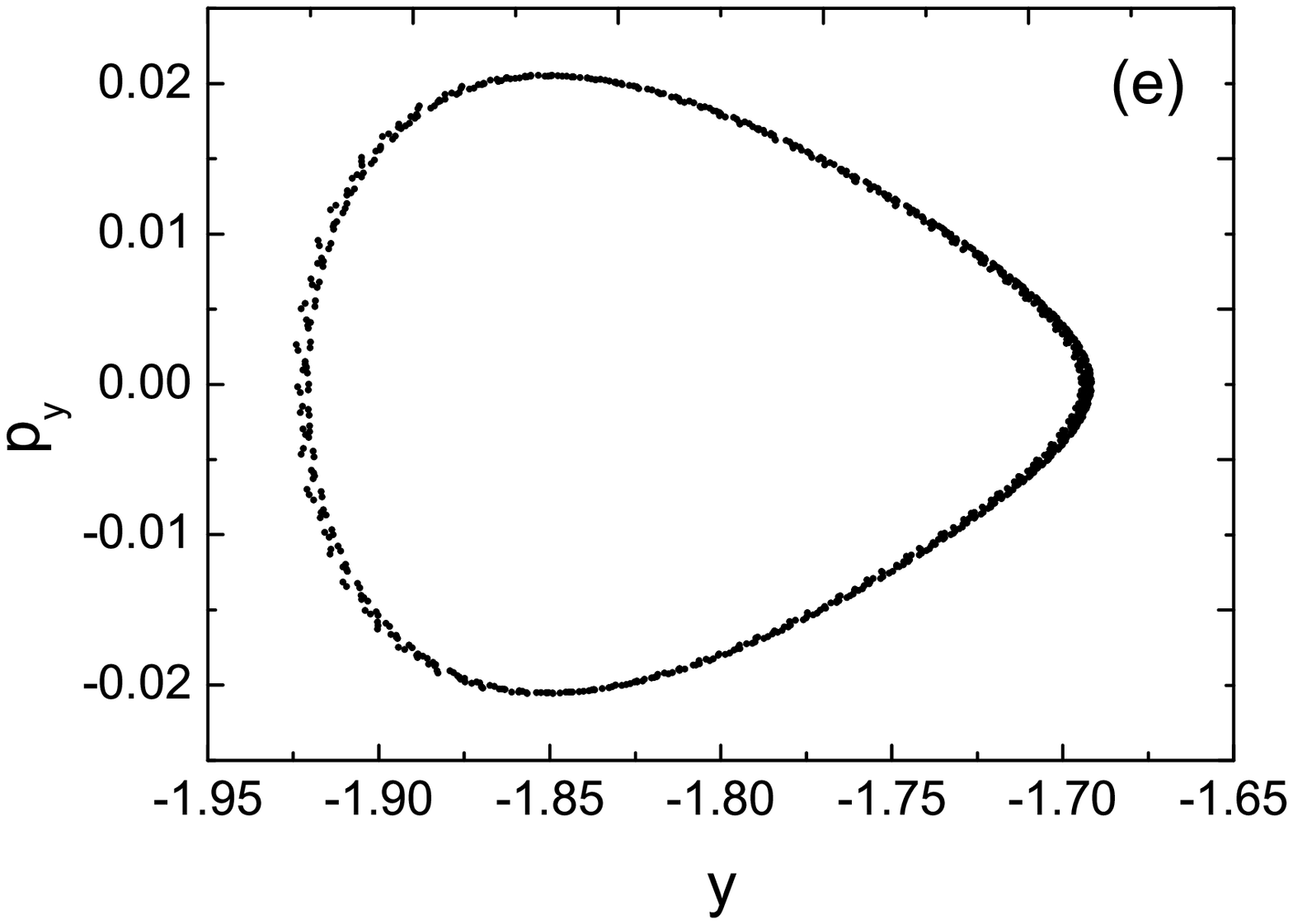}
\includegraphics[scale=0.26]{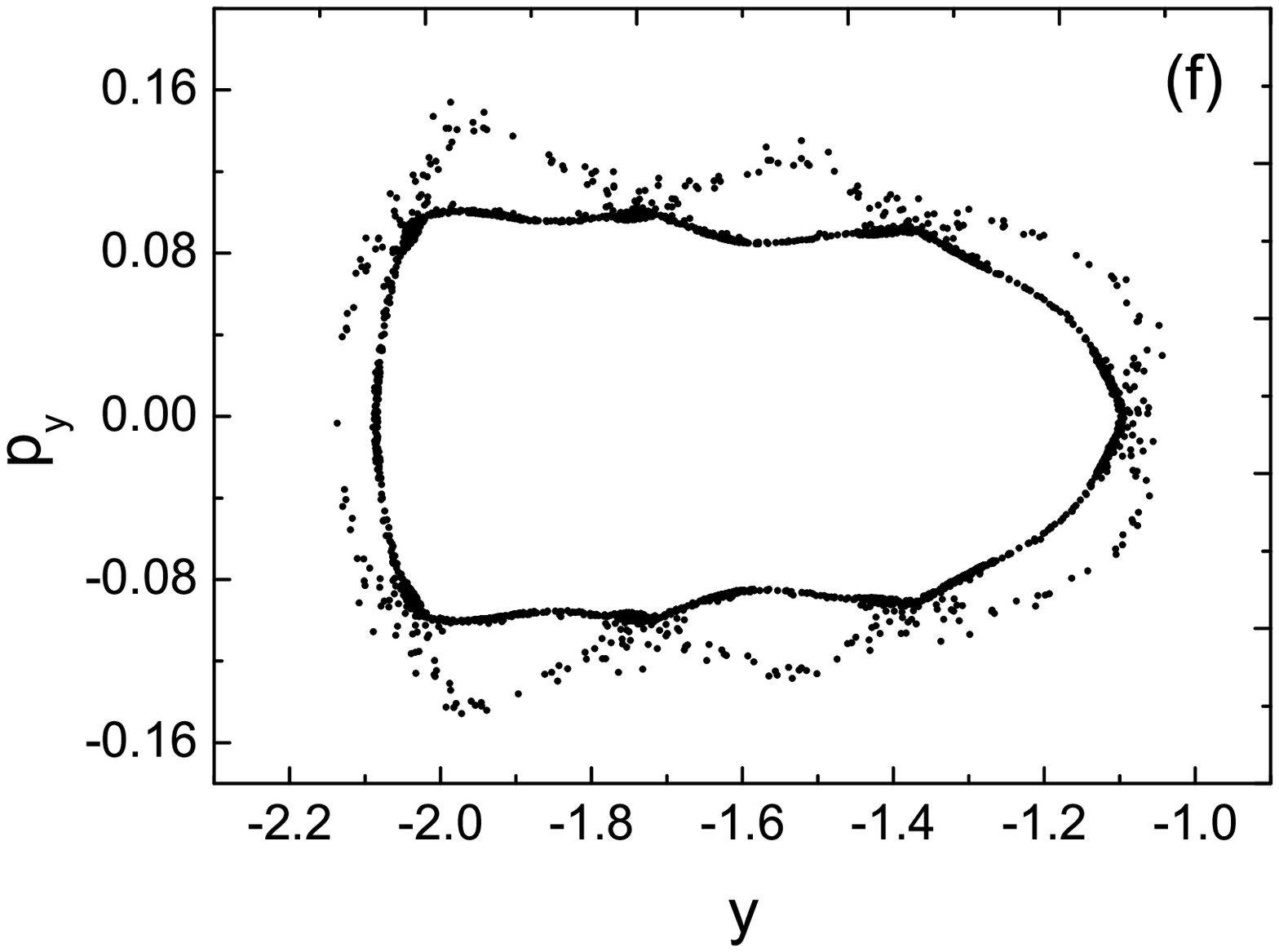}
\caption{Phase-space structures described by N4P with the larger
time step $\tau=0.1$ solving several orbits, which have only
different initial values $y$ with different starting values $p_x$.
(a) $y=-2.10$, (b) $y=-2.08$, (c) $y=-2.06$, (d) $y=-1.99$, (e)
$y=-1.6$, and (f) $y=-1.103$. There are regular many-islands tori
in (a)-(d). Only one single torus exists in (e). A chaotic
many-islands orbit is given in (f). All the dynamical structures
are physically performed by N4.  }} \label{f6}
\end{figure*}

\begin{figure*}
\center{
\includegraphics[scale=0.5]{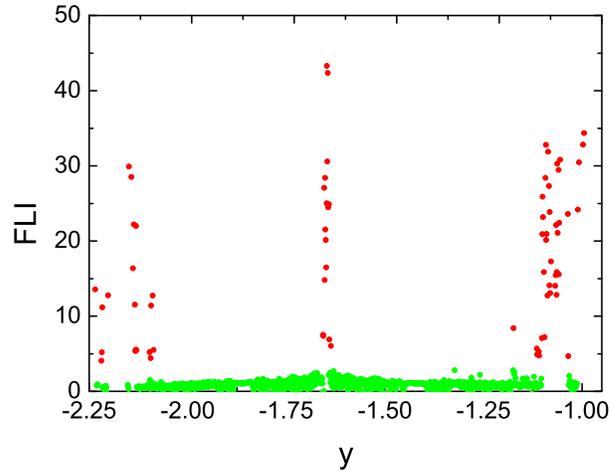}
\caption{Dependence of FLI on initial value $y$, described by
algorithm N4P with the larger time step $\tau=0.1$. Initial value
$x=0$ is fixed. Each of the FLIs is obtained after integration
time $t=3*10^{3}$. The FILs less than 4 correspond to the
regularity of orbits, whereas the FILs larger than 4 indicate the
chaoticity of orbits. Green corresponds to the regularity of
orbits, and Red indicates the chaoticity of orbits. Chaos mainly
occurs for the initial values of $y$ in the vicinity of
-2.25$\sim$-2.1, -1.6, and  -1.2$\sim$-1. }} \label{f7}
\end{figure*}

\begin{figure*}
\center{
\includegraphics[scale=0.35]{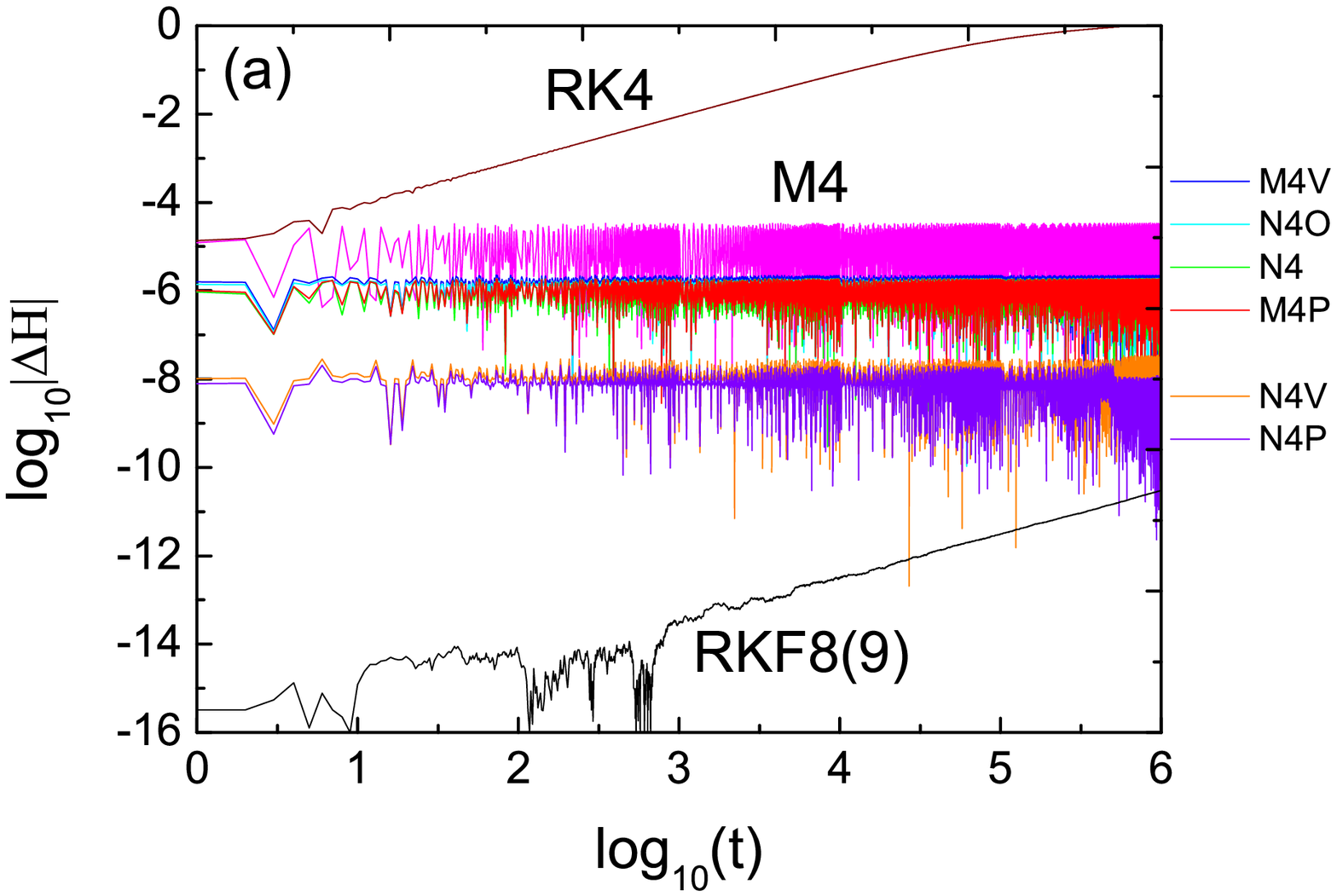}
\includegraphics[scale=0.35]{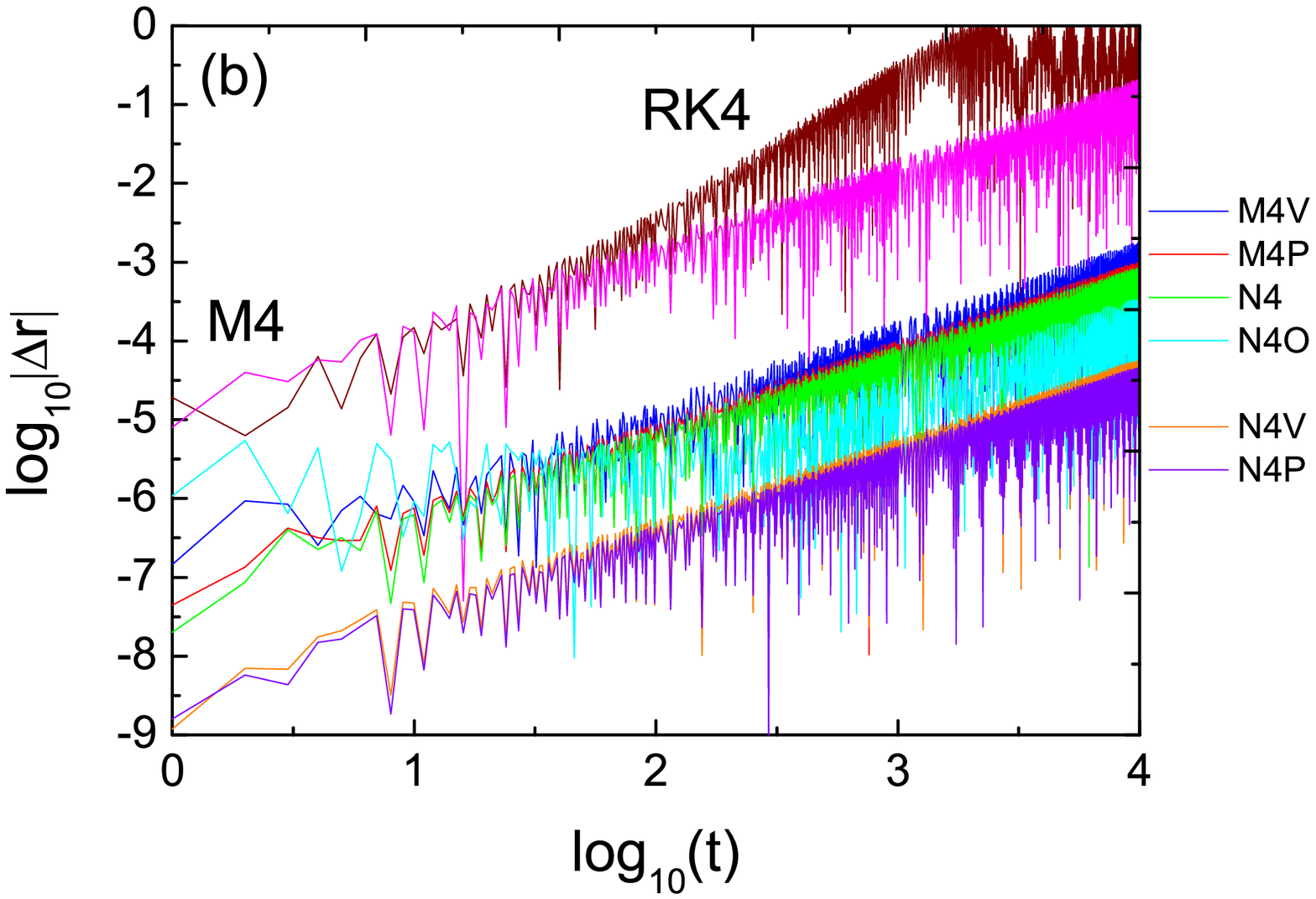}
\caption{Energy and position errors in  the spring pendulum
problem. The step size is $\tau=0.1$, and energy is $E=1/12$.
Initial conditions are $r=1.15$, $p_{r}=0$ and $\varphi=0.05\pi$.
(a) Energy errors for the algorithms. (b) Position errors for the
algorithms. These errors are clearly listed in Table 4. }}
\label{f8}
\end{figure*}

\begin{figure*}
\center{
\includegraphics[scale=0.26]{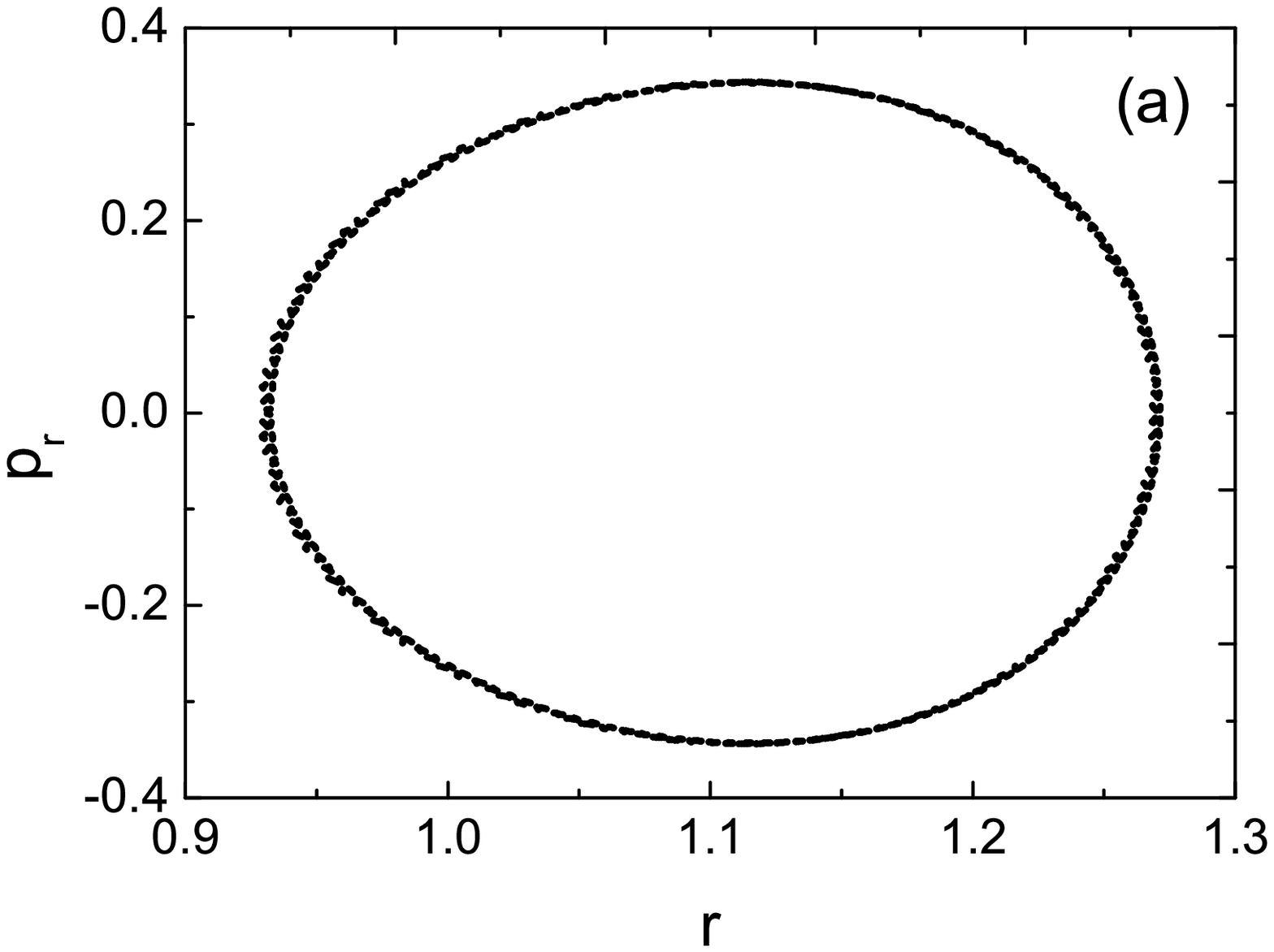}
\includegraphics[scale=0.26]{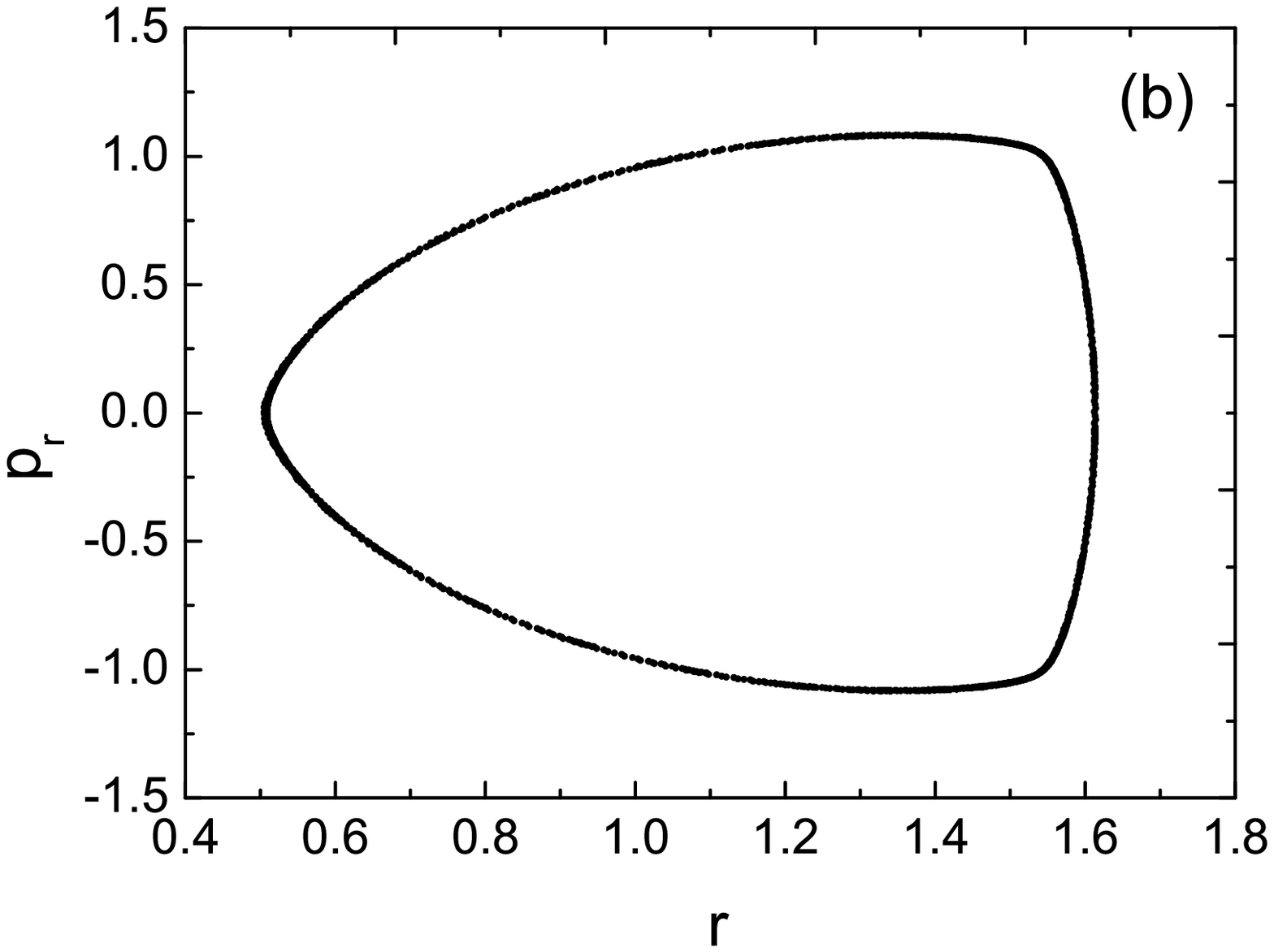}
\includegraphics[scale=0.26]{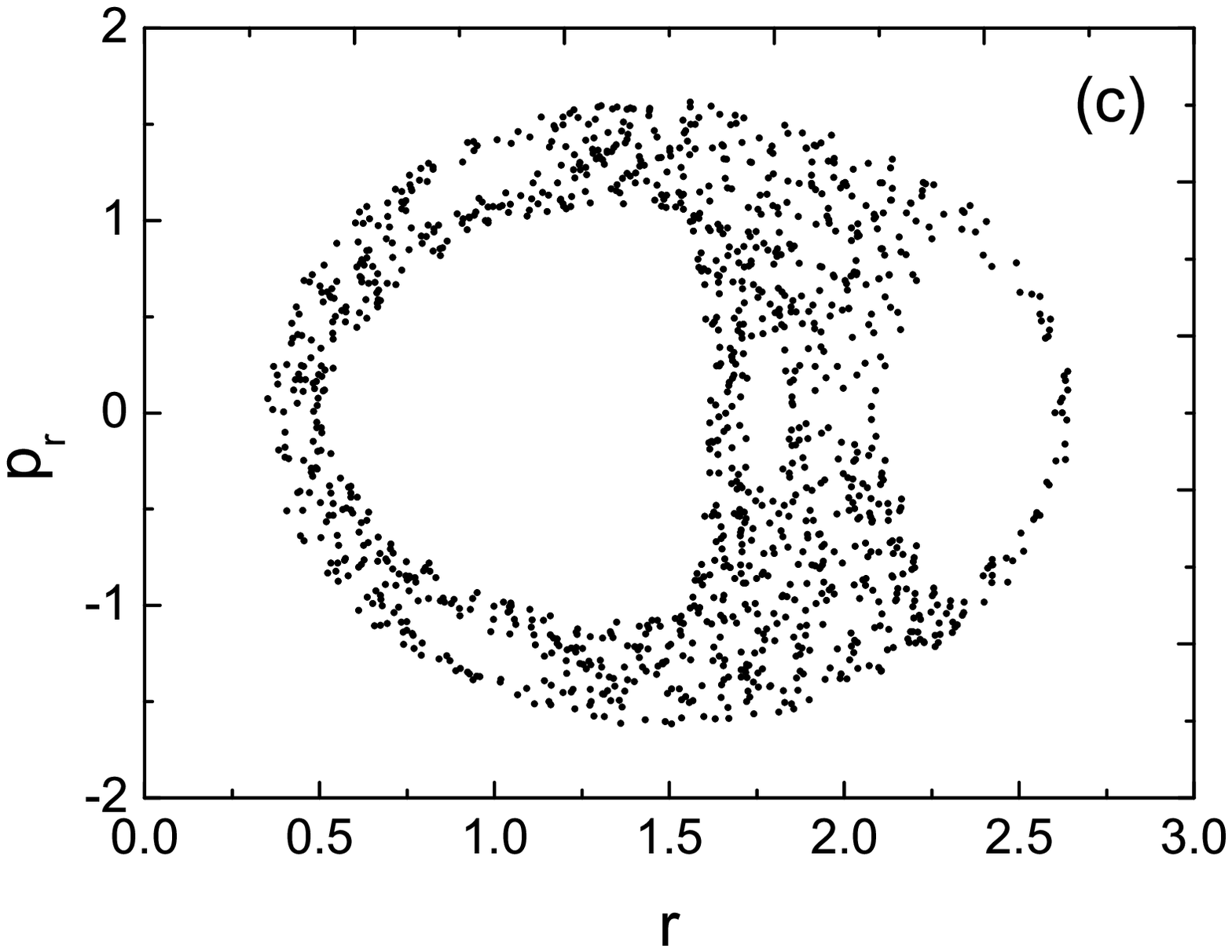}
\includegraphics[scale=0.26]{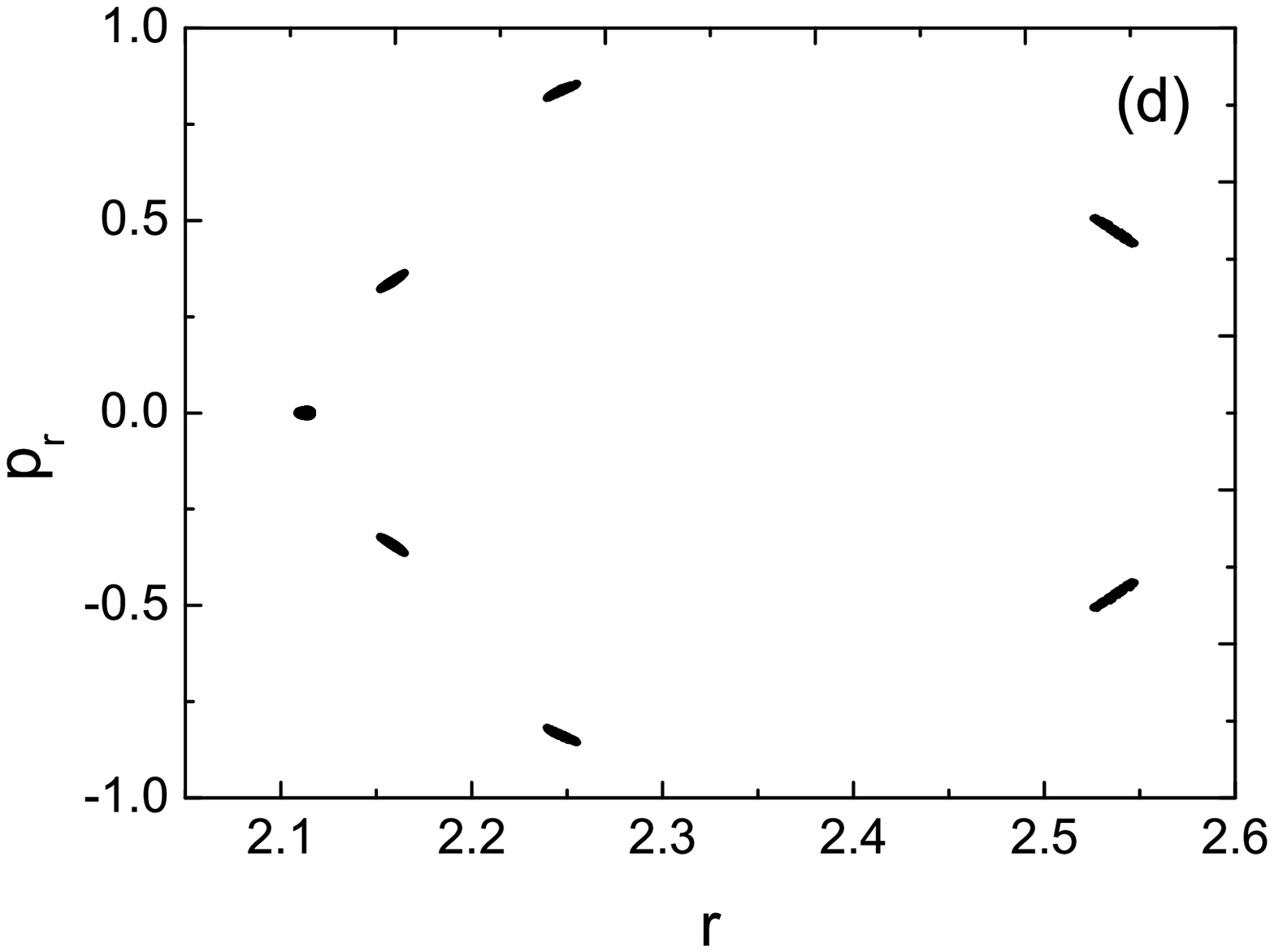}
\includegraphics[scale=0.26]{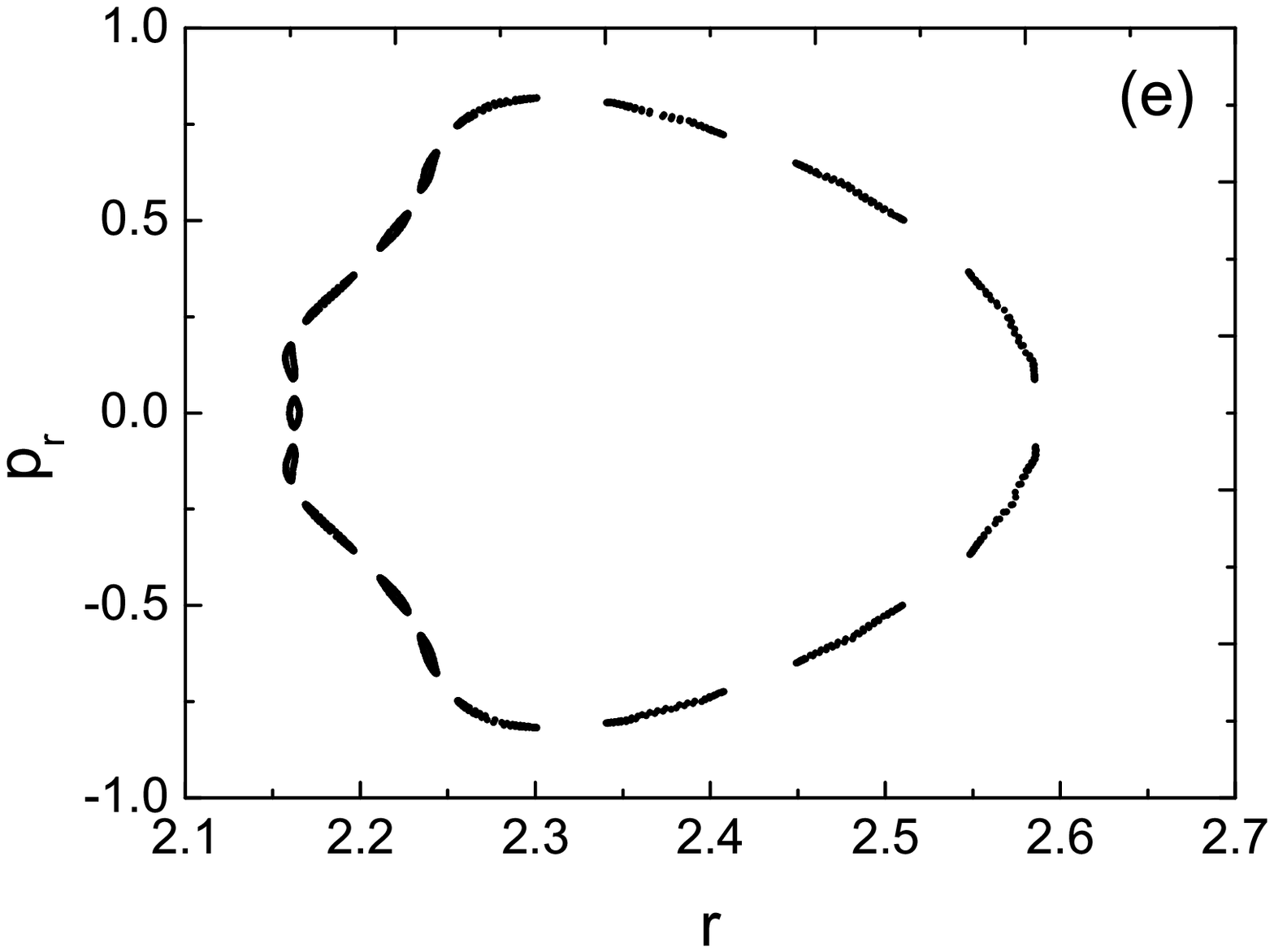}
\includegraphics[scale=0.26]{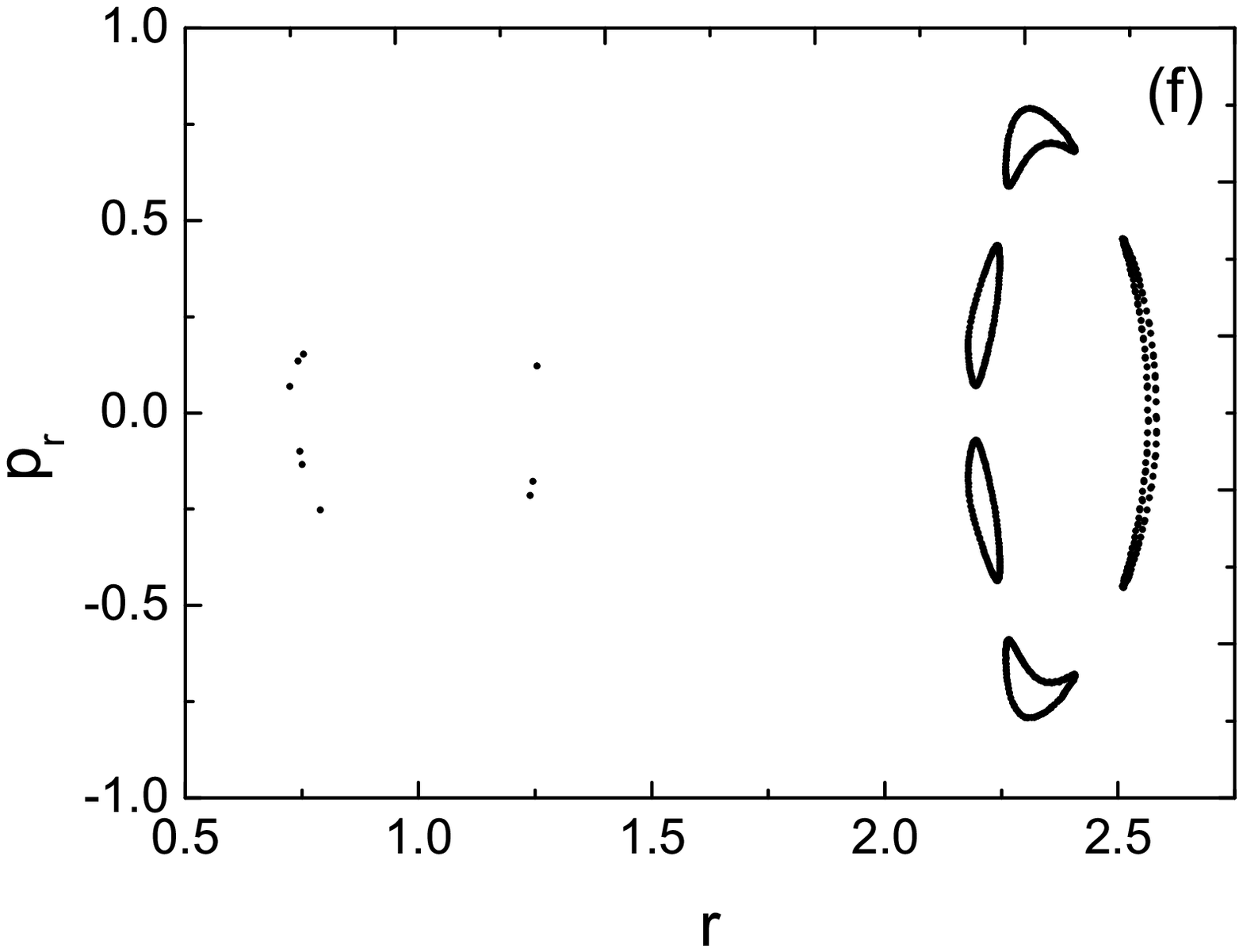}
\includegraphics[scale=0.26]{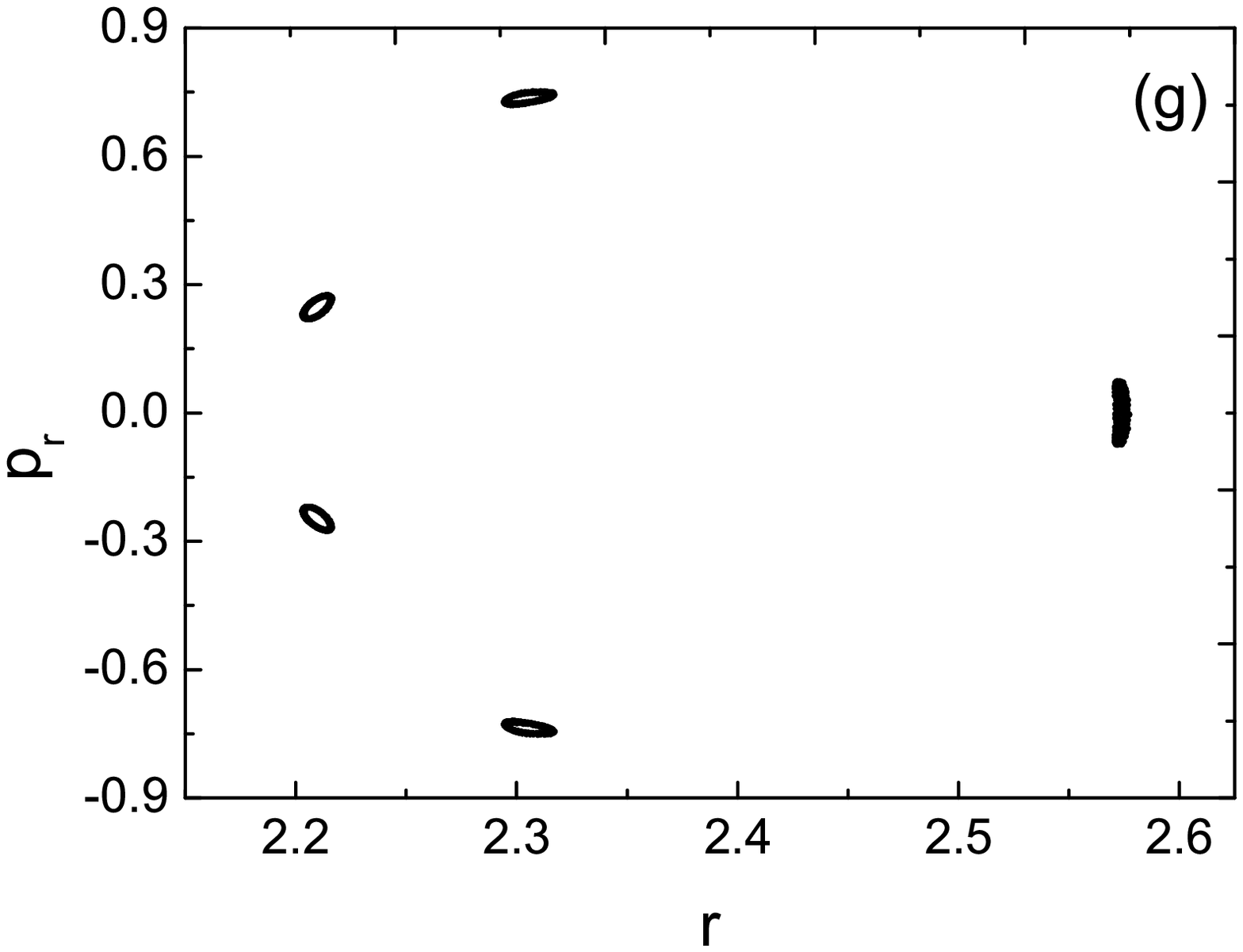}
\includegraphics[scale=0.26]{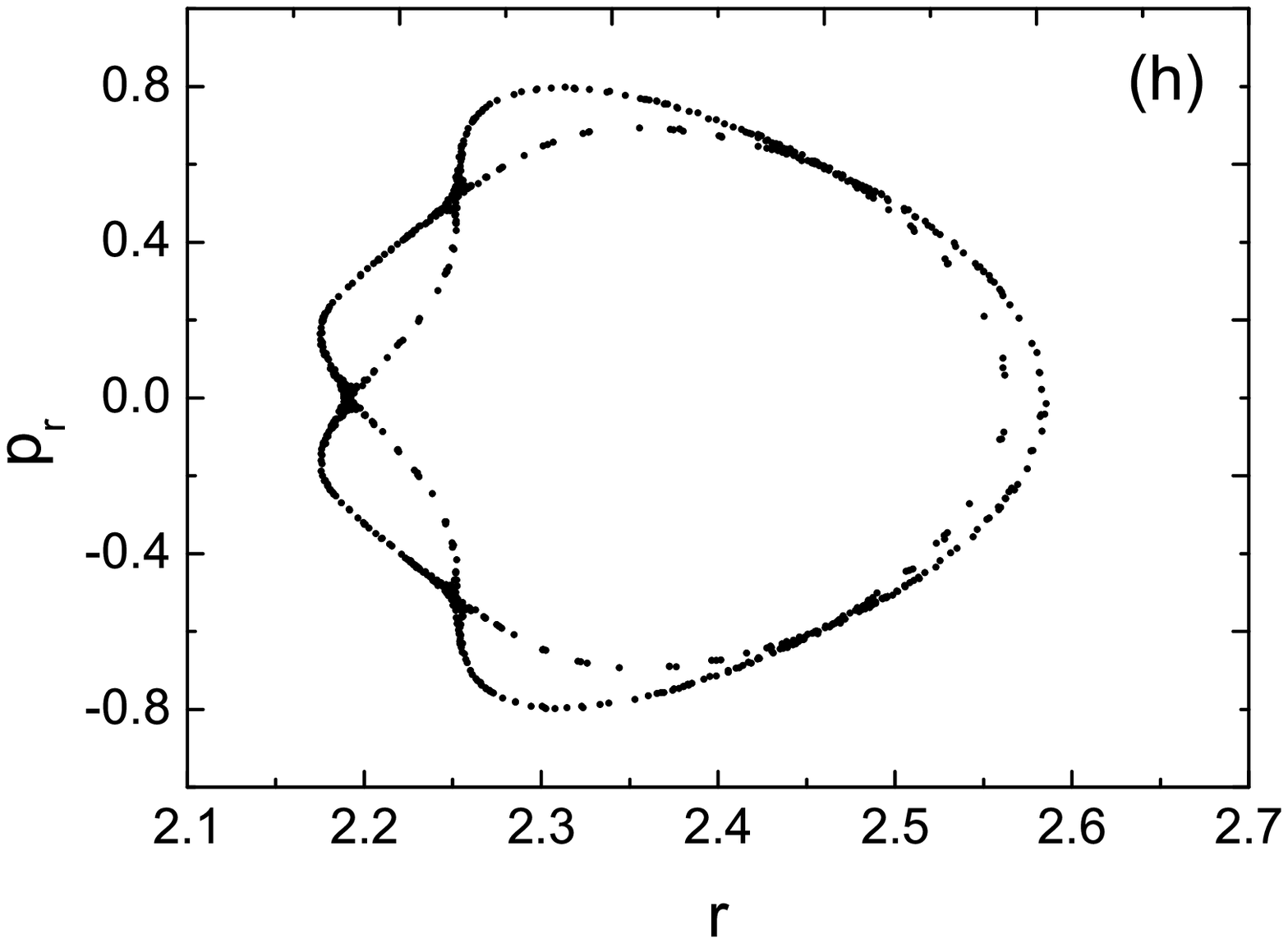}
\includegraphics[scale=0.26]{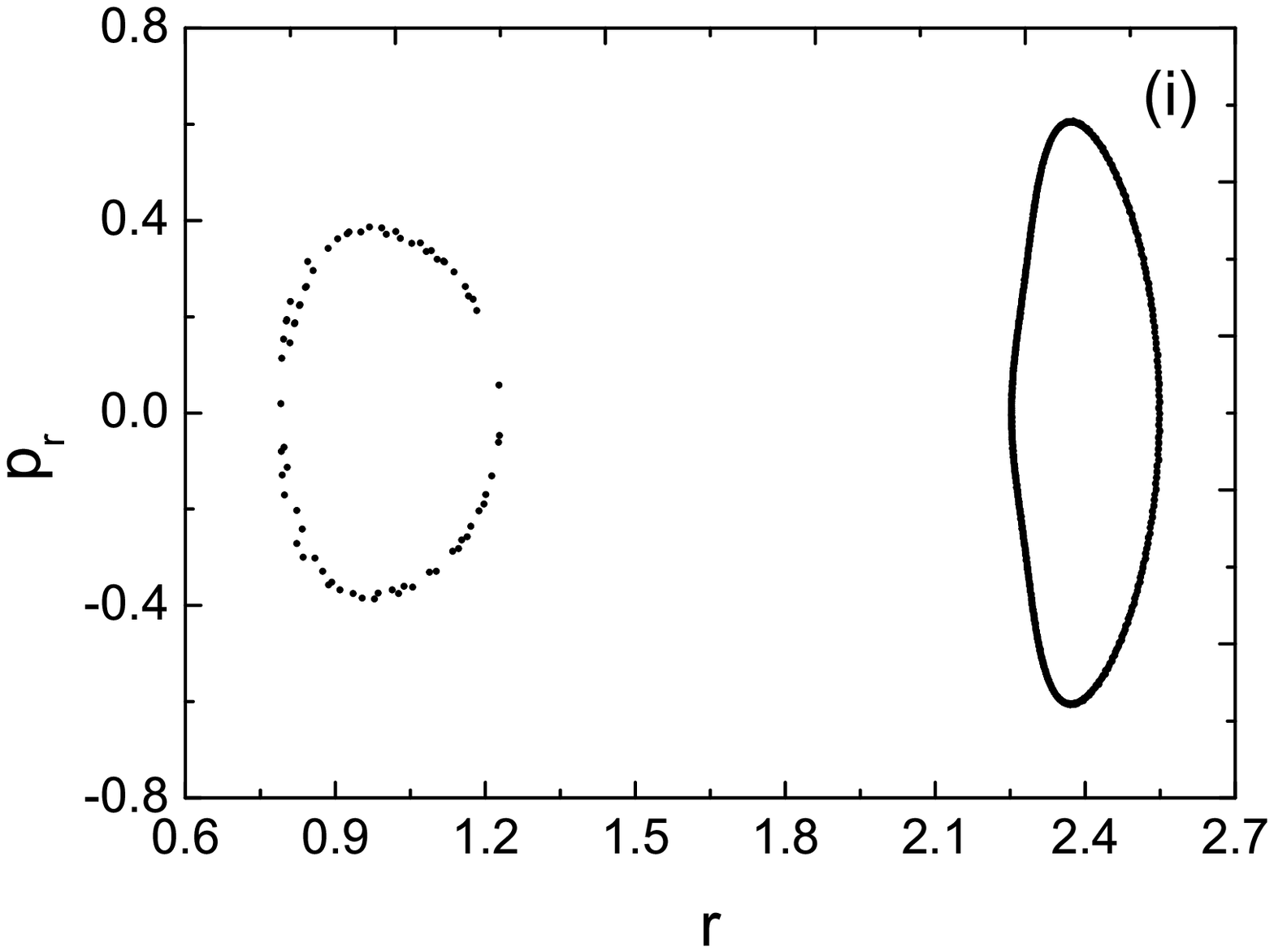}
\caption{Poincar\'{e} sections on the plane $\varphi=0$ and
$p_{\varphi}>0$ for various initial values of $\varphi$: (a)
$\varphi=0.05\pi$, (b) $\varphi=0.19\pi$, (c) $\varphi=0.2\pi$,
(d) $\varphi=0.35\pi$, (e) $\varphi=0.358\pi$, (f)
$\varphi=0.361\pi$, (g) $\varphi=0.366\pi$, (h)
$\varphi=0.376\pi$, and (i) $\varphi=0.39\pi$. The orbit in Fig. 7
is a regular single-torus orbit in Fig. 8(a). Regular single-torus
orbits also appear in (b) and (e). The orbits in (d), (f), (g) and
(i) are many-islands tori, but are chaotic in (c) and (h). These
results are consistently given by the methods M4, M4P, M4V, N4,
N4O, N4P, N4V and RKF8(9). }} \label{f9}
\end{figure*}

\begin{figure*}
\center{
\includegraphics[scale=0.35]{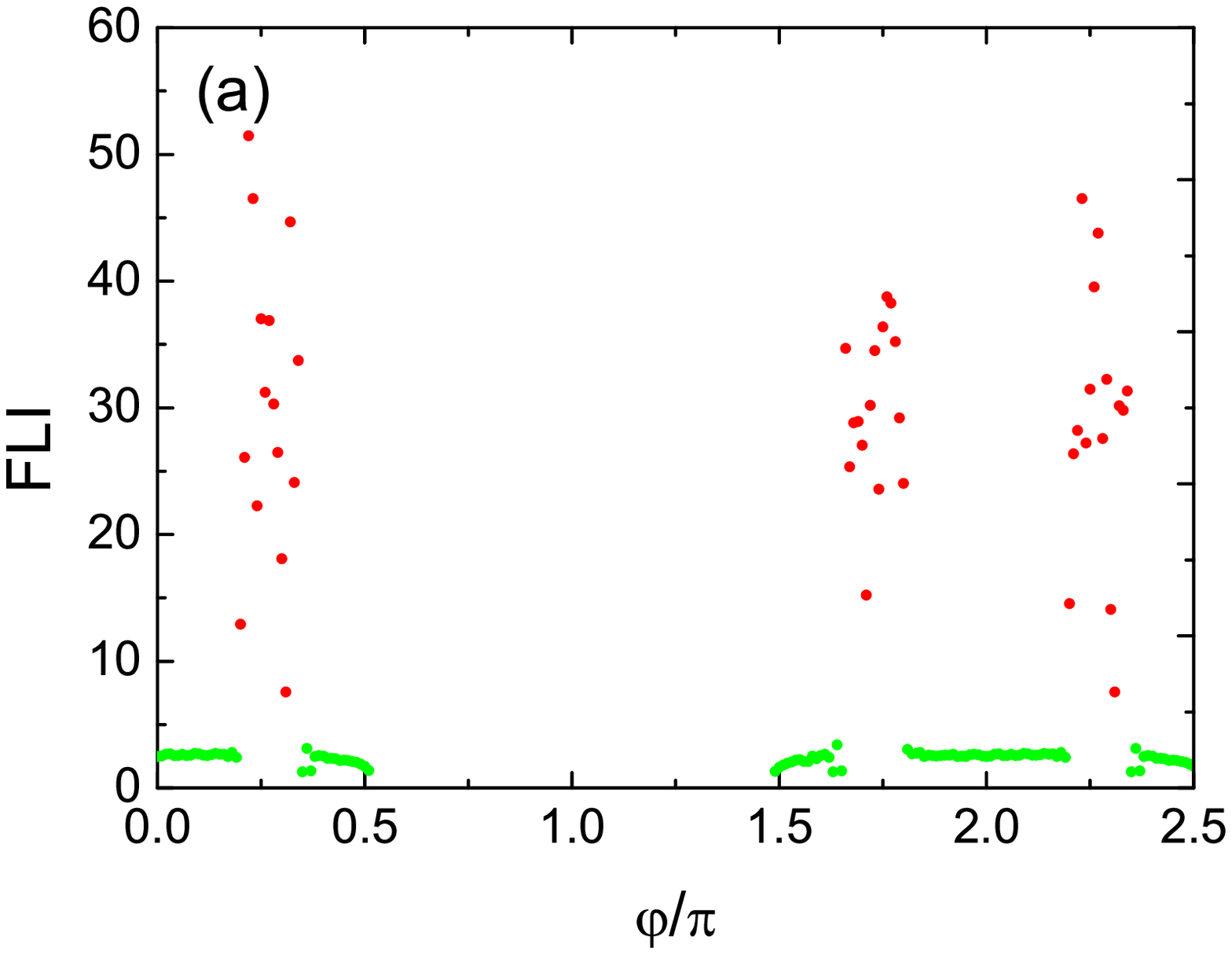}
\includegraphics[scale=0.35]{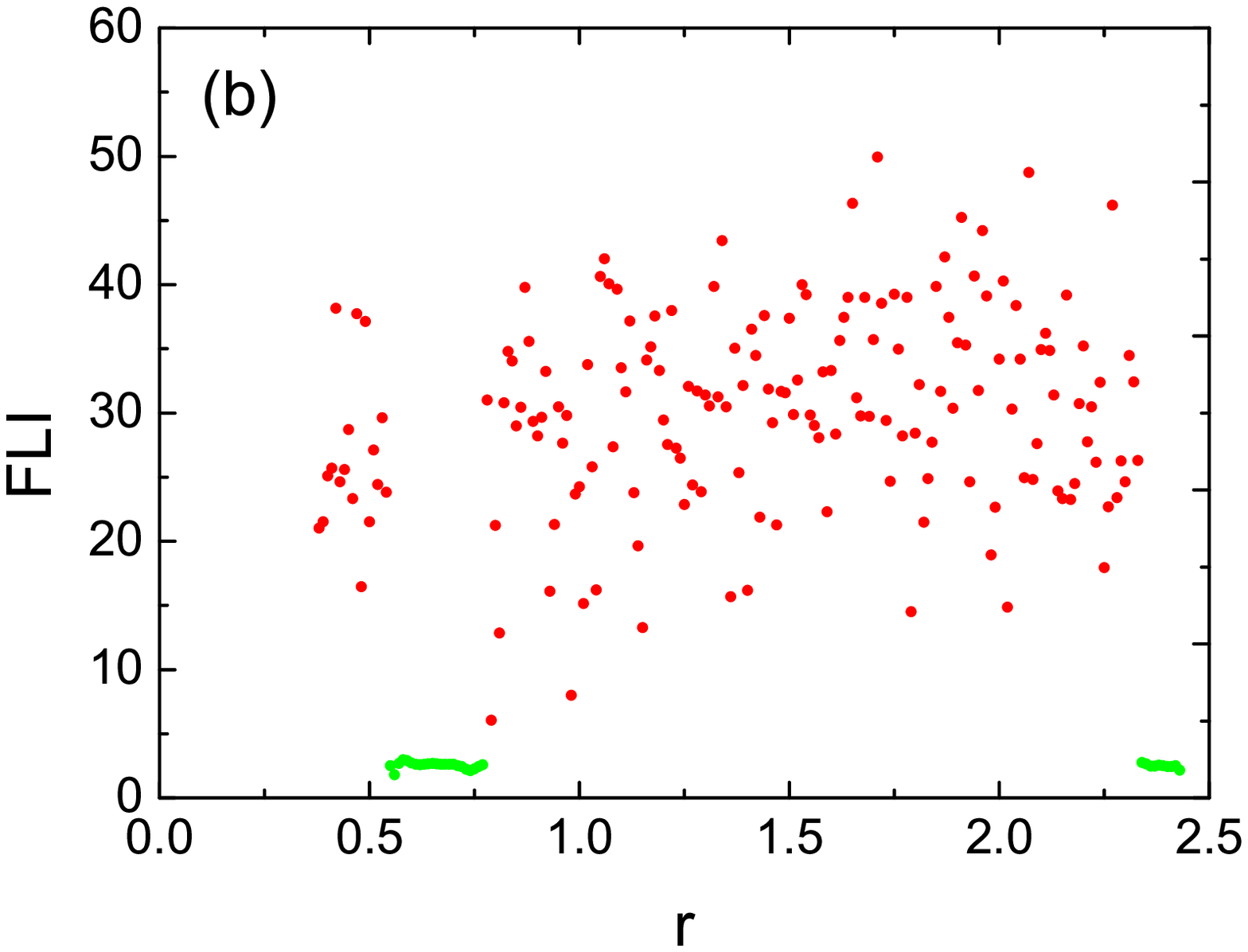}
\caption{Dependence of FLI on initial value $\varphi$ or $r$. (a)
Initial value $r=1.15$ is given. (b) Initial value
$\varphi=0.2\pi$ is fixed. Each of the FLIs is obtained by N4P
after integration time $t=1000$. The threshold of FLIs between the
ordered and chaotic cases is 4. Green corresponds to the
regularity of orbits, and Red indicates the chaoticity of orbits.
Chaos mainly occurs for the initial values of $\varphi$ in the
vicinity of 0.25, 1.75, and 2.25 in (a), and it does for the
initial values of $r$ in the vicinity of 0.5, and 0.75$\sim$2.25
in (b). }} \label{f10}
\end{figure*}

\begin{figure*}
\center{
\includegraphics[scale=0.35]{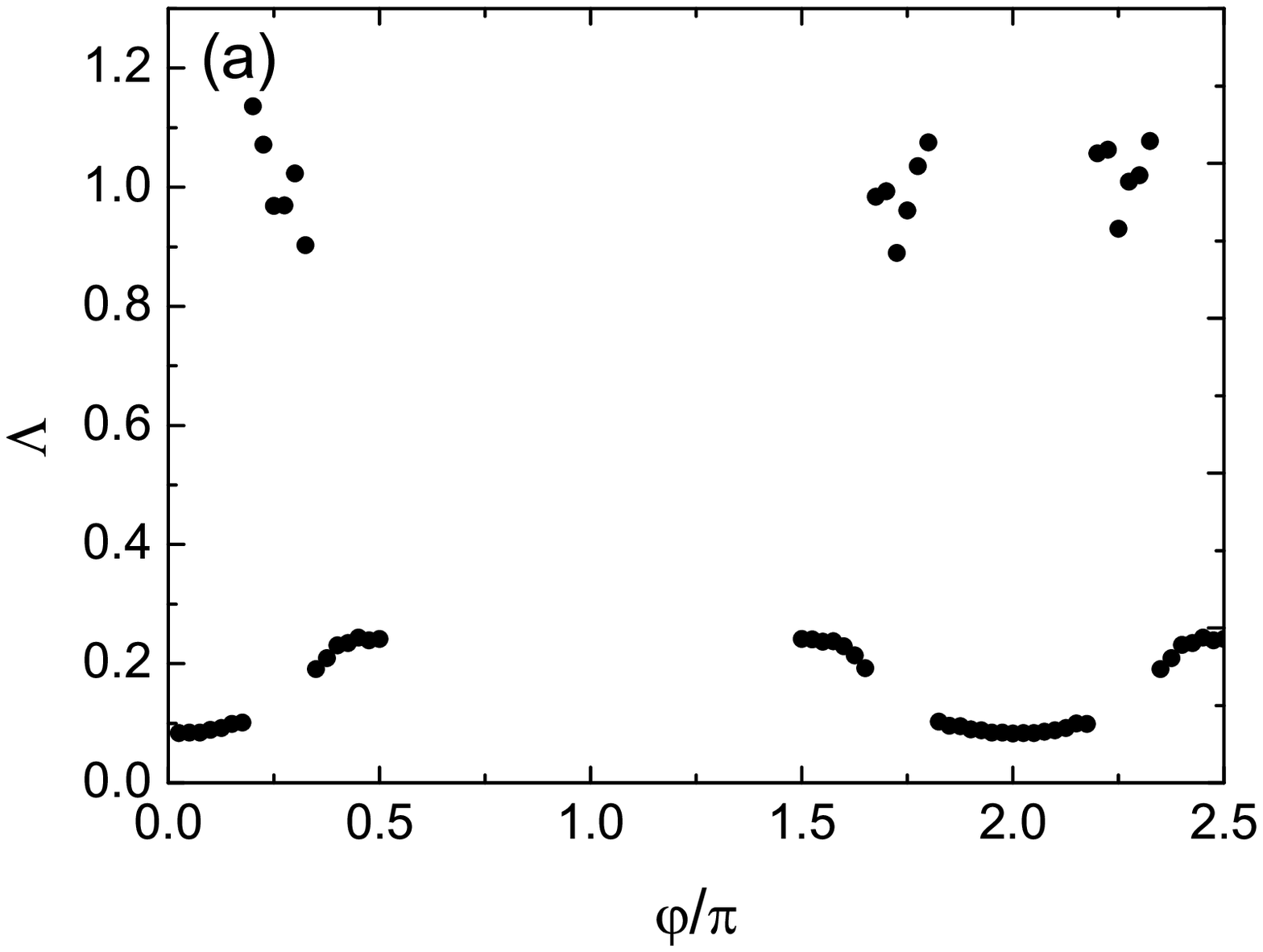}
\includegraphics[scale=0.35]{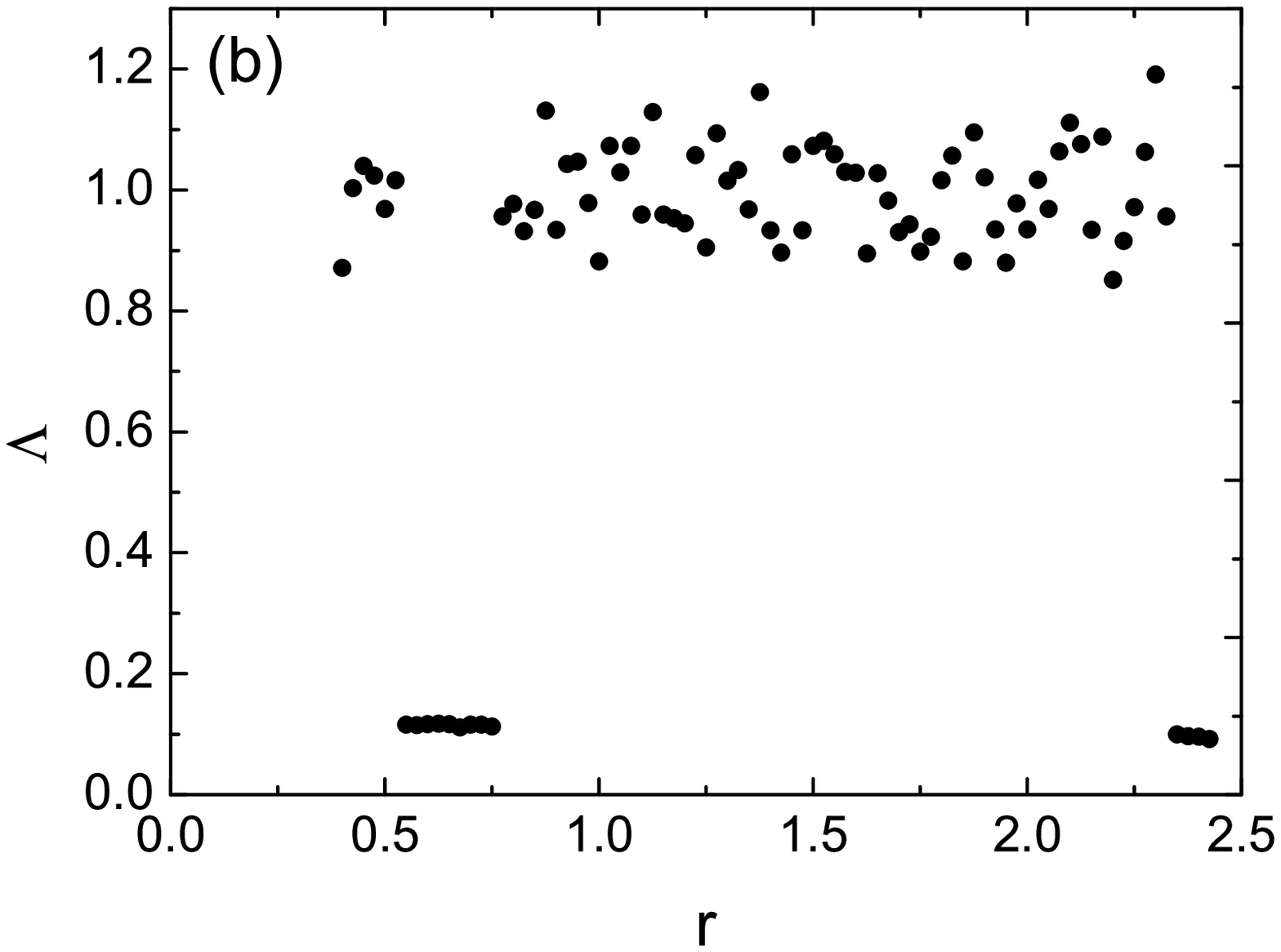}
\caption{ Same as Fig. 10 but the dependence of
the 0-1 test chaos indicator $\Lambda$ on  the initial value
$\varphi$ or $r$. The regular and chaotic properties described by
the 0-1 test are the same as those given by the FLIs in Fig. 10.
}} \label{f11}
\end{figure*}

%
%
%
%
%

\end{document}